\documentclass[a4paper,11pt]{article}

\usepackage{fullpage}
\usepackage{xurl}
\usepackage{graphicx}
\usepackage{alltt}
\usepackage{amsmath,amssymb,amsthm}
\usepackage[pdftex]{hyperref}
\usepackage[T1]{fontenc}
\usepackage[utf8]{inputenc}
\usepackage{lmodern}
\usepackage{csquotes}
\usepackage{cleveref}

\usepackage{amsfonts}
\usepackage{mathtools}
\usepackage{dsfont} 

\usepackage{enumitem} 

\usepackage[dvipsnames]{xcolor}
\usepackage[most]{tcolorbox}

\usepackage{color} 
\usepackage{soul} 

\usepackage{tikz}
\usepackage{pgfplots}
\usepackage{pgfplotstable}
\pgfplotsset{compat=1.18}

\usepackage{multirow} 
\usepackage{booktabs} 

\usepackage[english]{babel}

\newtheorem{thm}{Theorem}
\newtheorem{prp}{Proposition}
\newtheorem{lm}{Lemma}

\newtheorem{defi}{Definition}
\newtheorem{alg}{Algorithm}

\title{Triangle Factors in the Semi-Random Graph Process}
\author{Olaf Parczyk, Eva Schinzel}

\begin{document}

\maketitle

\begin{abstract}
The semi-random graph process is a single player game which starts with an empty graph on $n$ vertices and in each round the player is offered a random vertex and connects it to a vertex selected according to some strategy $\mathcal{S}$.
The goal is, in as few rounds as possible, to arrive at a graph that satisfies some graph property $\mathcal{P}$.
The  main contribution of this work is to show  that there is a strategy $\mathcal{S}$ which creates a $K_3$-factor in at most $2.929329 \, n$ rounds w.h.p.

\end{abstract}

\section{Introduction}\label{section_introduction}

The semi-random graph process (SRGP) is a 1-player game where, starting with the empty graph $G_0=([n],\emptyset)$, in each round the player is offered a uniformly random vertex $v_{t+1} \in [n]$ and then connects it to a vertex of choice $w_{t+1} \in [n]$, which creates $G_{t+1}:=G_t+v_{t+1}w_{t+1}$.
This model has strong connections to random graph models, positional games and, for example, the Achlioptas processes.
Since its introduction in 2018 in~\cite{ben2020semi} it has received a considerable amount of attention.
Typical questions in this model ask how long it takes such that the graph satisfies a certain graph property with high probability.
Most graph properties are various versions of connectivity, minimum degree, and the subgraph containment type (fixed size and spanning).
It turns out that a linear number of rounds is sufficient for these kinds of problems $\mathcal{P}$ and, hence, the goal is to determine the optimal limiting constant in front of $n$, denoted by $\tau_\mathcal{P}$.
Formally, for a monotone graph property $\mathcal{P}=\mathcal{P}(n)$ of $n$-vertex graphs we let
\begin{equation*}
  \tau_{\mathcal{P}} \coloneqq \lim_{q \to 1^-} \limsup_{n \to \infty}  \left( \frac{1}{n} \cdot \inf_{\mathcal{S}} \left( \min  \{t \geq 0:\mathbb{P}\left[ G_t^{\mathcal{S}}(n) \in \mathcal{P}(n) \right] \geq q \} \right) \right),
\end{equation*}
where $\mathcal{S}$ is any strategy\footnote{A strategy $\mathcal{S}$ is a sequence  $(f_{t+1})_{t=0}^{\infty}$ of distributions over $[n]$, s.t.\ in round $t$, where $t \in \mathbb{N}$, $w_{t}$ is chosen according to the distribution $f_{t}(v_1, w_1, ..., w_{t-1}, v_{t})$.} of the player and $G_t^{\mathcal{S}}(n)$ the $n$-vertex graph obtained in round $t$ when playing according to $\mathcal{S}$.

Focusing on spanning graphs, this has been intensively studied for perfect matchings $\mathcal{PM}$ and Hamilton cycles $\mathcal{HAM}$ with best bounds
\[0.93261 \le \tau_\mathcal{PM} \le 1.20524 \quad \text{and} \quad 1.26575 \le \tau_\mathcal{HAM} \le 1.81696\]
by different subsets of Frieze, Gao, MacRury, Pra\l{}at, and Sorkin~\cite{frieze2025hamilton,gao2022hamiltonfully,gao2022perfect}.
For any graph with maximum degree $\Delta$ an upper bound of $\Delta/2+o(\Delta)$ was shown~\cite{anastos2025optimally,ben2020very}, which requires $\Delta$ to be large and, hence, does not imply good bounds for triangle factors $\mathcal{K}_3^f$ or larger clique-factors $\mathcal{K}_k^f$.\footnote{A $K_k$-factor $\mathcal{K}_k^f$ is the disjoint union of $\lfloor n/k \rfloor$ copies of $K_k$ and at most $k-1$ isolated vertices.}

The main result of this paper is an improved upper bound for $\mathcal{K}_3^f$ generalizing the case $\tau_{\mathcal{K}_2}^f=\tau_{\mathcal{PM}}$.
\begin{thm}\label{main_theorem_k3_factor_upper}
$\tau_{\mathcal{K}_3^f} \leq 2.929329$.
\end{thm}
This result also appears in the master thesis of the second author~\cite{eva_thesis}, where also further generalizations to $\tau_{\mathcal{K}_k^f}$ are discussed and, in particular, $\tau_{\mathcal{K}_4^f} \leq 4.890484$.\footnote{Note that actually in~\cite{eva_thesis} the bounds are slightly weaker. The new bounds are solely due to an improved implementation.}
We give the proof of \Cref{main_theorem_k3_factor_upper} in \Cref{k3_factor_section}.
The analytical and numerical solutions of differential equations involved in the proof were computed using SymPy (cf.~\cite{sympyreference}) and SciPy (cf.~\cite{scipyreference}).
The corresponding Python code is provided on \href{https://github.com/Mysterious-Rodent-Princess/clique-factors}{github}.

As for lower bounds on $\tau_{\mathcal{K}_3^f}$, fortunately, the proof of the lower bound on $\tau_{\mathcal{HAM}}$ from \cite{gao2022hamiltonfully} only uses that the graph is $2$-regular.
Hence, it also applies to triangle factors and $\tau_{\mathcal{K}_3^f} \ge 1.26575$ follows. In these types of problems, a triangle factor is harder to find than a Hamilton cycle and hence it would be interesting to improve that lower bound.

The SRGP also has an offline version, i.e.~the player knows the entire sequence $v_1, v_2, v_3, ...$ of random vertices before making the choices.
The quantity $\tau'_{\mathcal{P}}$ is the analog of $\tau_{\mathcal{P}}$ for the offline setting.
In the offline setting $\tau'_{\mathcal{P}}$ only depends on the number of offerings of each vertex and not the order.
\begin{thm}\label{main_theorem_k3_factor_offline}
$\tau'_{\mathcal{K}_3^f}=x \approx 1.1462$, where $x$ is the positive solution to $e^{-x}(x + 2)=1$.
\end{thm}

The proof essentially boils down to the observation that there are two ways to build triangles, one vertex offered twice and one vertex once or all vertices offered once, and a short optimization.
We do not give the proof here and instead refer to~\cite{eva_thesis}, where it is also extended to larger cliques.
We note that this quantity also shows up in the proof of the lower bound on $\tau_{\mathcal{HAM}}$ in \cite[Theorem 1.2]{frieze2025hamilton}.

\section{Proof overview and tools}
\label{section_balls}

To give an upper bound of the form $\tau_{\mathcal{P}} \le x$ one needs to provide a strategy $\mathcal{S}$ and prove that with high probability (w.h.p.) within $xn$ rounds $G_t^S$ satisfies $\mathcal{P}$.
A naive strategy to build a perfect matching would be to maintain a matching $M$ and if a vertex $v \notin V(M)$ is offered, to connect it to another vertex $w\notin V(M)$, add $vw$ to $M$, and pick an arbitrary $w$ if $v \in M$.
This way it takes $\Theta(n \log n)$ until $M$ is a perfect matching.
Hence, any strategy that constructs a perfect matching in linear time has to use some form of augmentation, i.e.~larger flexible structures that might change the current matching.
For example, if the two endpoints of a matching edge are offered one can create a path of length $3$, which contains a larger matching that incorporates two new vertices, none of which had to be offered.
Doing this in parallel and in different phases is how the strategy for the bound $\tau_\mathcal{PM} \le 1.20524 $ from~\cite{gao2022perfect} is achieved.

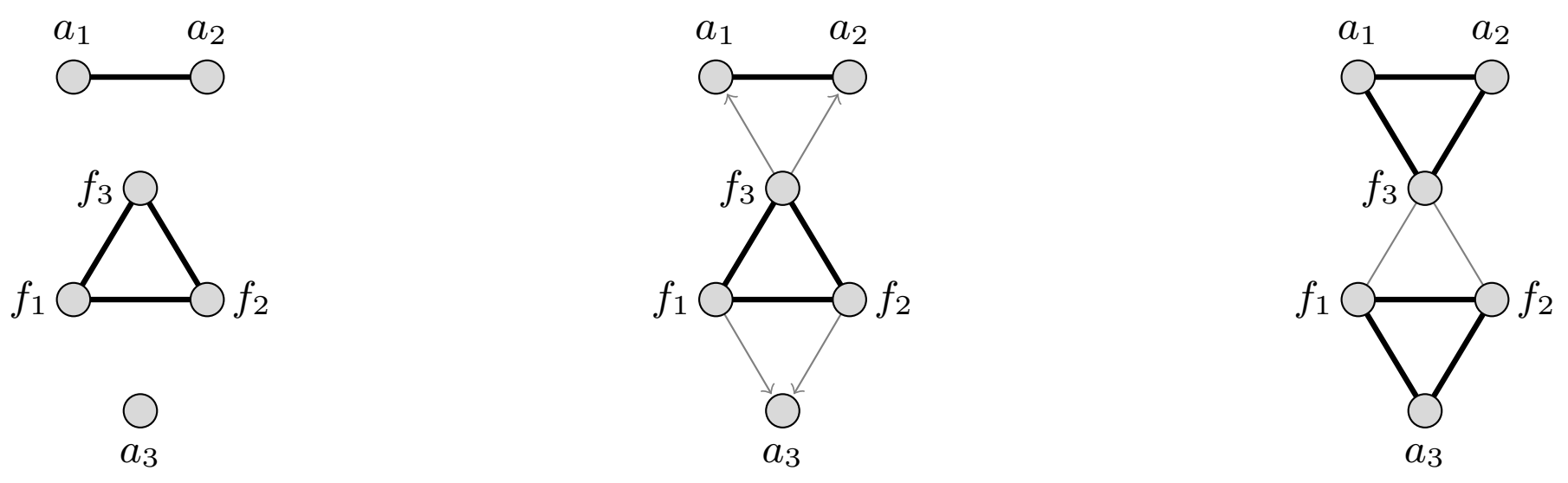
\begin{figure}[htbp]
\begin{center}
 \begin{minipage}{0.3\textwidth}
 \begin{center}
\begin{tikzpicture}[scale=0.85]
\node at (-0.6,3.4) {\footnotesize $a_1$};
\node at (0.6,3.4) {\footnotesize $a_2$};
\node at (0,-0.4) {\footnotesize $a_3$};
\node at (-0.4,2) {\footnotesize $f_3$};
\node at (-1,1) {\footnotesize $f_1$};
\node at (1,1) {\footnotesize $f_2$};
\draw[very thick] (-0.6,1) -- (0,2);
\draw[very thick] (0.6,1) -- (0,2);
\draw[very thick] (-0.6,1) -- (0.6,1);
\draw[very thick] (-0.6,3) -- (0.6,3);
\draw[fill=gray!30] (0,0) circle (1.5mm);
\draw[fill=gray!30] (-0.6,1) circle (1.5mm);
\draw[fill=gray!30] (0.6,1) circle (1.5mm);
\draw[fill=gray!30] (0,2) circle (1.5mm);
\draw[fill=gray!30] (-0.6,3) circle (1.5mm);
\draw[fill=gray!30] (0.6,3) circle (1.5mm);
\end{tikzpicture}
\end{center}
\end{minipage}
 \begin{minipage}{0.3\textwidth}
 \begin{center}
\begin{tikzpicture}[scale=0.85]
\node at (-0.6,3.4) {\footnotesize $a_1$};
\node at (0.6,3.4) {\footnotesize $a_2$};
\node at (0,-0.4) {\footnotesize $a_3$};
\node at (-0.4,2) {\footnotesize $f_3$};
\node at (-1,1) {\footnotesize $f_1$};
\node at (1,1) {\footnotesize $f_2$};
\draw[gray, ->] (-0.6,1) -- (-0.1,0.15);
\draw[gray, ->] (0.6,1) -- (0.1,0.15);
\draw[very thick] (-0.6,1) -- (0,2);
\draw[very thick] (0.6,1) -- (0,2);
\draw[very thick] (-0.6,1) -- (0.6,1);
\draw[gray, ->] (0,2) -- (-0.5,2.85);
\draw[gray, ->] (0,2) -- (0.5,2.85);
\draw[very thick] (-0.6,3) -- (0.6,3);
\draw[fill=gray!30] (0,0) circle (1.5mm);
\draw[fill=gray!30] (-0.6,1) circle (1.5mm);
\draw[fill=gray!30] (0.6,1) circle (1.5mm);
\draw[fill=gray!30] (0,2) circle (1.5mm);
\draw[fill=gray!30] (-0.6,3) circle (1.5mm);
\draw[fill=gray!30] (0.6,3) circle (1.5mm);
\end{tikzpicture}
\end{center}
\end{minipage}
\begin{minipage}{0.3\textwidth}
\begin{center}
\begin{tikzpicture}[scale=0.85]
\node at (-0.6,3.4) {\footnotesize $a_1$};
\node at (0.6,3.4) {\footnotesize $a_2$};
\node at (0,-0.4) {\footnotesize $a_3$};
\node at (-0.4,2) {\footnotesize $f_3$};
\node at (-1,1) {\footnotesize $f_1$};
\node at (1,1) {\footnotesize $f_2$};
\draw[very thick] (0,0) -- (-0.6,1);
\draw[very thick] (0,0) -- (0.6,1);
\draw[gray] (-0.6,1) -- (0,2);
\draw[gray] (0.6,1) -- (0,2);
\draw[very thick] (-0.6,1) -- (0.6,1);
\draw[very thick] (0,2) -- (-0.6,3);
\draw[very thick] (0,2) -- (0.6,3);
\draw[very thick] (-0.6,3) -- (0.6,3);
\draw[fill=gray!30] (0,0) circle (1.5mm);
\draw[fill=gray!30] (-0.6,1) circle (1.5mm);
\draw[fill=gray!30] (0.6,1) circle (1.5mm);
\draw[fill=gray!30] (0,2) circle (1.5mm);
\draw[fill=gray!30] (-0.6,3) circle (1.5mm);
\draw[fill=gray!30] (0.6,3) circle (1.5mm);
\end{tikzpicture}
\end{center}
\end{minipage}
\end{center}
    \caption{On the left $A$ and $F$ before additional offerings, in the middle four suitable offerings and the respective choices with directed edges, and on the right the successful augmentation.}
    \label{fig:augmentation_triangles}
\end{figure}

Now we set out to generalize this concept of an augmenting path to triangles.
Suppose $A$ is an induced subgraph on three vertices of $G_t$ for some step $t$  with at least one and at most two edges, say $V(A) = \{a_1, a_2, a_3\}$, and $\{a_1, a_2\} \in E(A)$.
Moreover, let $F$ be a triangle in $G_t$ with $V(F) = \{f_1, f_2, f_3\}$.
Now suppose that $f_1$ and $f_2$ are offered once each and both connected to $a_3$.
Suppose in addition that $f_3$ is offered twice, once connected to $a_1$ and once to $a_2$.
Then there are two triangles on vertex sets $\{a_1,a_2,f_3\}$ and $\{a_3,f_1,f_2\}$.
This way, using four semi-random edges, two new triangles can be used instead of the old one, cf.~\Cref{fig:augmentation_triangles}, and we call this an \emph{augmentation}.
Note that we will never use the second edge even if it is there, but because we will need to accurately track their numbers we distinguish between \emph{$A^1$-structures}, which have only the edge $a_1a_2$ induced in $\{a_1,a_2,a_3\}$ and \emph{$A^2$-structures} which have $a_1a_2$ and one other edge.
We call both of these \emph{augmentable structures} and, in particular, in this context we will abuse notation and often only refer to the edge-set as a representative of the whole structure.
Naturally, the first vertex of a triangle $F$ that is offered gets the role of $f_3$, because it needs another offering, and the roles of $f_1$ and $f_2$ are assigned to the other two vertices of $F$ arbitrarily.
Similarly, when $f_3$ is offered we call $a_1$ the \emph{appropriate vertex} for $f_3$ and analogously afterwards the \emph{appropriate vertices} for $f_3$, $f_2$, and $f_1$ are $a_2$, $a_3$, and $a_3$, respectively.
Once a vertex of $F$ was offered and connected to an augmentable structure $A$, we say that $F$ is \emph{responsible} for $A$.
To make this approach efficient we will later have many triangles that are responsible for the same $A$.

\subsection{Strategy}\label{strategy_description}

Now we are ready to describe the strategy for our proof of \Cref{main_theorem_k3_factor_upper}.
We note that it is sufficient to do this for $n$ that is divisible by three, so from now on we assume that this is the case.
Our goal is to quickly build triangles and use them for augmentation.
We will begin with a simple Phase 1:
\begin{enumerate}[label=(\roman*)]
    \item Whenever an unsaturated vertex is offered, we connect it to another unsaturated vertex to get a matching edge.
    \item When an endpoint of a matching edge is offered, we connect it to an unsaturated vertex and get an $A^2$-structure.
    \item When an endpoint of an $A^2$-structure is offered, we create a triangle.
    \item Otherwise we make an arbitrary move that we forget.
\end{enumerate}

\begin{figure}[htbp]
    \begin{center}
 \begin{tikzpicture}[scale=0.8]
   \draw[thick](-8.25, 0.75) -- (-8.25,1.25) -- (-7.75,1.25) -- (-7.75,0.75) -- cycle;
\draw[fill=gray!30] (-8,1) circle (1.5mm);
\draw[fill=gray!30] (-8.5, 0.2) circle (1.5mm);
 \draw[very thick](-5.5, 0.2) -- (-5,1);
\draw[fill=gray!30] (-5,1) circle (1.5mm);
\draw[fill=gray!30] (-5.5, 0.2) circle (1.5mm);
   \draw[thick](-1.75, -0.05) -- (-1.75,0.45) -- (-1.25,0.45) -- (-1.25,-0.05) -- cycle;
 \draw[very thick](-1.5, 0.2) -- (-2,1);
\draw[fill=gray!30] (-2,1) circle (1.5mm);
\draw[fill=gray!30] (-1,1) circle (1.5mm);
\draw[fill=gray!30] (-1.5, 0.2) circle (1.5mm);
  \draw[very thick](1.5, 0.2) -- (1,1);
    \draw[very thick](1.5, 0.2) -- (2,1);
\draw[fill=gray!30] (1,1) circle (1.5mm);
\draw[fill=gray!30] (2,1) circle (1.5mm);
\draw[fill=gray!30] (1.5, 0.2) circle (1.5mm);
    \draw[thick](4.75, 0.75) -- (4.75,1.25) -- (5.25,1.25) -- (5.25,0.75) -- cycle;
 \draw[very thick](5.5, 0.2) -- (6,1);
  \draw[very thick](5.5, 0.2) -- (5,1);
\draw[fill=gray!30] (5,1) circle (1.5mm);
\draw[fill=gray!30] (6,1) circle (1.5mm);
\draw[fill=gray!30] (5.5, 0.2) circle (1.5mm);
 \draw[very thick](8,1) -- (9,1);
 \draw[very thick](8.5, 0.2) -- (8,1);
  \draw[very thick](8.5, 0.2) -- (9,1);
\draw[fill=gray!30] (8,1) circle (1.5mm);
\draw[fill=gray!30] (9,1) circle (1.5mm);
\draw[fill=gray!30] (8.5, 0.2) circle (1.5mm);
   \node at (-6.8,0.8) {\footnotesize \textit{(i)}};
\draw[gray, ->] (-7.3,0.5) -- (-6.3,0.5);
   \node at (0,0.8) {\footnotesize \textit{(ii)}};
\draw[gray, ->] (-0.5,0.5) -- (0.5,0.5);
   \node at (7,0.8) {\footnotesize \textit{(iii)}};
\draw[gray, ->] (6.5,0.5) -- (7.5,0.5);
  \end{tikzpicture}
 \end{center}
    \caption{(i) turns two unsaturated vertices into a matching edge, (ii) turns a matching edge and one unsaturated vertex into an $A^2$-structure, and (iii) turns an $A^2$-structure into a triangle. The squares indicate the offered vertices.}
    \label{fig:phase_1_simple}
\end{figure}
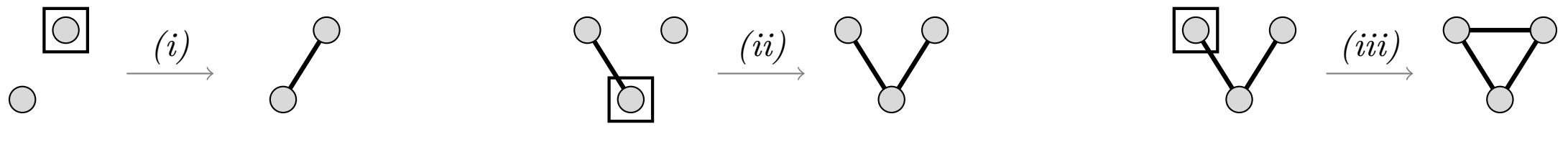

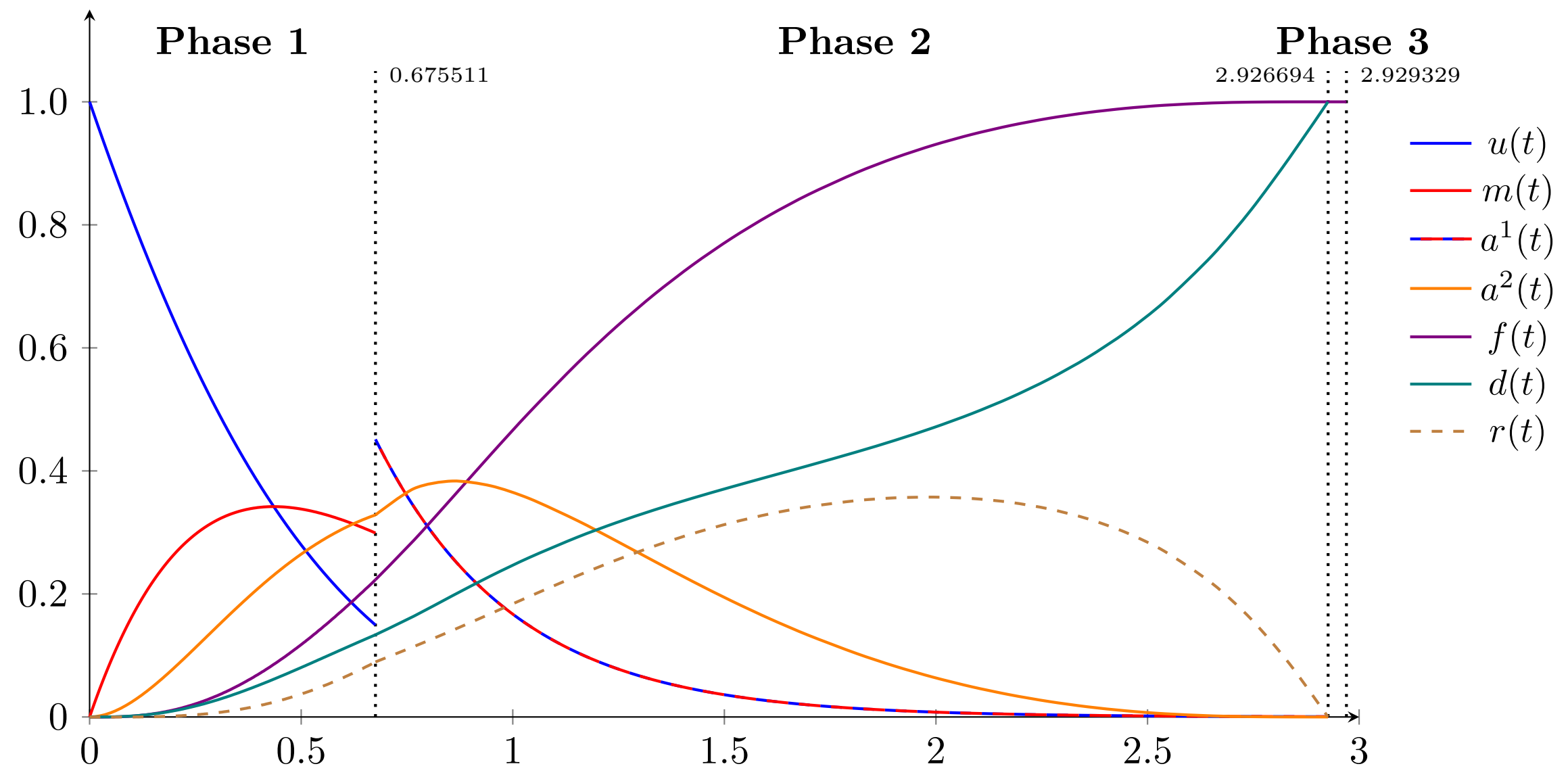
\begin{figure}[htbp]
    \centering
 \begin{tikzpicture}

\begin{axis}[
    width=14cm,
    height=8.5cm,
    xmin=0, xmax=3,
    ymin=0, ymax=1.15,
    xlabel={},
    ylabel={},
    axis lines=left,
    xtick={0,0.5,1,1.5,2,2.5,3},
    ytick={0,0.2,0.4,0.6,0.8,1.0},
    yticklabels={0,0.2,0.4,0.6,0.8,1.0},
    tick style={thin},
    clip=false,
    legend pos=outer north east,
    legend style={font=\small, draw=none, yshift=-1cm},
    every axis plot/.append style={line width=0.8pt},
]

\def\xI{0.675511}
\def\xII{2.926694}   
\def\xIII{2.970}     

\addplot[dotted, thick, black, forget plot]
    coordinates {(\xI,0) (\xI,1.05)};
\addplot[dotted, thick, black, forget plot]
    coordinates {(\xII,0) (\xII,1.05)};
\addplot[dotted, thick, black, forget plot]
    coordinates {(\xIII,0) (\xIII,1.05)};

\node[above, font=\tiny, rotate=0, anchor=south west] at (axis cs:\xI, 1.01)
    {$0.675511$};
\node[above, font=\tiny, rotate=0, anchor=south west] at (axis cs:\xII-0.3, 1.01)
    {$2.926694$};
\node[above, font=\tiny, rotate=0, anchor=south west] at (axis cs:\xIII, 1.01)
    {$2.929329$};

\node[font=\bfseries] at (axis cs:0.338, 1.10) {Phase 1};
\node[font=\bfseries] at (axis cs:1.808, 1.10) {Phase 2};
\node[font=\bfseries] at (axis cs:2.985, 1.10) {Phase 3};


\addplot[domain=0:\xI, samples=200, color=blue]
    { cos(deg(sqrt(2)*x)) * exp(-2*x) };
\addlegendentry{$u(t)$}

\addplot[domain=0:\xI, samples=200, color=red]
    { sqrt(2) * sin(deg(sqrt(2)*x)) * exp(-2*x) };
\addlegendentry{$m(t)$}

\addplot[domain=\xI:\xII, samples=200, color=blue, dashed, dash pattern=on 6pt off 3pt]
    { 0.45136 * exp(-3.06232*(x - \xI)) };
\addlegendentry{$a^1(t)$}
\addplot[domain=\xI:\xII, samples=200, color=red, dashed,  dash pattern=on 6pt off 3pt, dash phase=6pt, forget plot]
    { 0.45136 * exp(-3.06232*(x - \xI)) };

\addplot[domain=0:\xI, samples=200, color=orange]
    { (3 - 3*cos(deg(sqrt(2)*x))) * exp(-2*x) };
\addlegendentry{$a^2(t)$}

\addplot[domain=0:\xI, samples=200, color=violet]
    { (2*cos(deg(sqrt(2)*x)) - sqrt(2)*sin(deg(sqrt(2)*x)) - 3) * exp(-2*x) + 1 };
\addlegendentry{$f(t)$}

\addplot[domain=0:\xI, samples=200, color=teal]
    { 0.25
      - (14/17*cos(deg(sqrt(2)*x)) + 2*sqrt(2)/17*sin(deg(sqrt(2)*x))) * exp(-2*x)
      + (39/68*cos(deg(sqrt(15)/2*x)) + 27*sqrt(15)/340*sin(deg(sqrt(15)/2*x))) * exp(-3.5*x) };
\addlegendentry{$d(t)$}

\addplot[domain=0:\xI, samples=200, color=brown, dashed]
    { 3/8
      + (6/17*cos(deg(sqrt(2)*x)) - 21*sqrt(2)/17*sin(deg(sqrt(2)*x))) * exp(-2*x)
      - (99/136*cos(deg(sqrt(15)/2*x)) - 57*sqrt(15)/680*sin(deg(sqrt(15)/2*x))) * exp(-3.5*x) +
      3/8
      + (42/17*cos(deg(sqrt(2)*x)) + 6*sqrt(2)/17*sin(deg(sqrt(2)*x)) - 3) * exp(-2*x)
      + (21/136*cos(deg(sqrt(15)/2*x)) - 111*sqrt(15)/680*sin(deg(sqrt(15)/2*x))) * exp(-3.5*x)};
\addlegendentry{$r(t)$}


\addplot[color=violet, smooth] coordinates {
  (0.67551, 0.22307)
  (0.76645, 0.28732)
  (0.85740, 0.35674)
  (0.94834, 0.42700)
  (1.03928, 0.49530)
  (1.15296, 0.57530)
  (1.24391, 0.63389)
  (1.33485, 0.68723)
  (1.42579, 0.73526)
  (1.51674, 0.77809)
  (1.60768, 0.81600)
  (1.69862, 0.84928)
  (1.81230, 0.88492)
  (1.90325, 0.90906)
  (1.99419, 0.92968)
  (2.08513, 0.94707)
  (2.17608, 0.96151)
  (2.26702, 0.97324)
  (2.35796, 0.98250)
  (2.44891, 0.98950)
  (2.53985, 0.99445)
  (2.65353, 0.99816)
  (2.74447, 0.99953)
  (2.83542, 0.99996)
  (2.92636, 0.99999)
  (2.970, 1.00000)
};

\addplot[color=orange, smooth] coordinates {
  (0.67551, 0.32837)
  (0.76645, 0.37126)
  (0.85740, 0.38355)
  (0.94834, 0.37568)
  (1.03928, 0.35505)
  (1.15296, 0.31912)
  (1.24391, 0.28650)
  (1.33485, 0.25298)
  (1.42579, 0.22006)
  (1.51674, 0.18872)
  (1.60768, 0.15953)
  (1.69862, 0.13285)
  (1.81230, 0.10322)
  (1.90325, 0.08254)
  (1.99419, 0.06449)
  (2.08513, 0.04898)
  (2.17608, 0.03591)
  (2.26702, 0.02514)
  (2.35796, 0.01654)
  (2.44891, 0.00999)
  (2.53985, 0.00530)
  (2.65353, 0.00177)
  (2.74447, 0.00046)
  (2.83542, 0.00004)
  (2.92636, 0.00000)
};

\addplot[color=teal, smooth] coordinates {
  (0.67551, 0.13361)
  (0.76645, 0.16390)
  (0.85740, 0.19707)
  (0.94834, 0.22941)
  (1.03928, 0.25925)
  (1.15296, 0.29235)
  (1.24391, 0.31567)
  (1.33485, 0.33667)
  (1.42579, 0.35589)
  (1.51674, 0.37393)
  (1.60768, 0.39130)
  (1.69862, 0.40857)
  (1.81230, 0.43080)
  (1.90325, 0.44974)
  (1.99419, 0.47035)
  (2.08513, 0.49325)
  (2.17608, 0.51923)
  (2.26702, 0.54921)
  (2.35796, 0.58417)
  (2.44891, 0.62554)
  (2.53985, 0.67493)
  (2.65353, 0.75064)
  (2.74447, 0.82421)
  (2.83542, 0.90906)
  (2.92636, 0.99998)
};

\addplot[color=brown, dashed, smooth] coordinates {
(0.67551, 0.08946)
(0.76645, 0.11449)
(0.85740, 0.14074)
(0.94834, 0.16814)
(1.03928, 0.19590)
(1.15296, 0.22964)
(1.24391, 0.25488)
(1.33485, 0.27785)
(1.42579, 0.29822)
(1.51674, 0.31570)
(1.60768, 0.33022)
(1.69862, 0.34169)
(1.81230, 0.35172)
(1.90325, 0.35621)
(1.99419, 0.35741)
(2.08513, 0.35511)
(2.17608, 0.34890)
(2.26702, 0.33820)
(2.35796, 0.32230)
(2.44891, 0.29993)
(2.53985, 0.26944)
(2.65353, 0.21660)
(2.74447, 0.15933)
(2.83542, 0.08649)
(2.92636, 0.00002)
};
\end{axis}

\draw[red, thick] (13.02,4.675) -- (13.236,4.675);
\draw[red, thick] (13.33,4.675) -- (13.52,4.675);

\end{tikzpicture}
    \caption{The functions represent in step $tn$ the fraction of vertices that are unsaturated $u(t)$, in matching edges $m(t)$, in triangles $f(t)$, or in an augmentable structure $a^1(t)$ resp.~$a^2(t)$. 
    At the end of Phase 1 we have $m(t)=2u(t)$ and they are merged into $a^1(t)$.
    Additionally, we distinguish between the vertices in triangles that have a responsibility, $r(t)$, and those which do not, $d(t)$.}
    \label{fig:plot}
\end{figure}

We will do this until the number of isolated vertices equals the number of matching edges, because then they exactly form $A^1$-structures and this happens w.h.p.\ at some point when three divides $n$.
This first phase will roughly take $0.675511n$ rounds and at this point we will have roughly $0.074358 n$ triangles.
In \Cref{fig:plot} the development of the fraction of vertices in the different roles over the duration of Phase 1 is plotted.

In Phase 2, which includes augmentation, we proceed as follows:
\begin{enumerate}[label=(\roman*)]
    \item Whenever an endpoint of an $A^1$-structure is offered, we create an $A^2$-structure.
    \item When an endpoint of an $A^2$-structure is offered, we create a triangle.
    \item When a vertex in a triangle without responsibility is offered, we randomly choose an augmentable structure and start the augmentation as described above.
    \item[(iv)-(vi)] When a suitable vertex in a triangle with responsibility is offered, we continue the augmentation procedure.
    \item[(vii)] Otherwise we make an arbitrary move that we forget.
\end{enumerate}

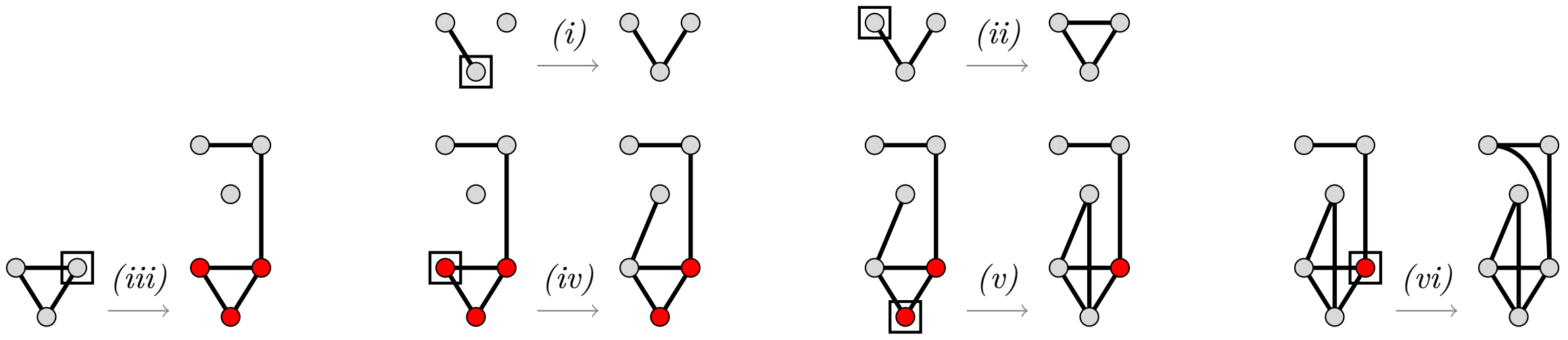
\begin{figure}[htbp]
    \begin{center}
 \begin{tikzpicture}[scale=0.6]
     \node at (0,0.8) {\footnotesize \textit{(i)}};
\draw[gray, ->] (-0.5,0.3) -- (0.5,0.3);
   \draw[thick](-1.75, -0.05) -- (-1.75,0.45) -- (-1.25,0.45) -- (-1.25,-0.05) -- cycle;
 \draw[very thick](-1.5, 0.2) -- (-2,1);
\draw[fill=gray!30] (-2,1) circle (1.5mm);
\draw[fill=gray!30] (-1,1) circle (1.5mm);
\draw[fill=gray!30] (-1.5, 0.2) circle (1.5mm);
  \draw[very thick](1.5, 0.2) -- (1,1);
    \draw[very thick](1.5, 0.2) -- (2,1);
\draw[fill=gray!30] (1,1) circle (1.5mm);
\draw[fill=gray!30] (2,1) circle (1.5mm);
\draw[fill=gray!30] (1.5, 0.2) circle (1.5mm);
     \node at (7,0.8) {\footnotesize \textit{(ii)}};
\draw[gray, ->] (6.5,0.3) -- (7.5,0.3);
    \draw[thick](4.75, 0.75) -- (4.75,1.25) -- (5.25,1.25) -- (5.25,0.75) -- cycle;
 \draw[very thick](5.5, 0.2) -- (6,1);
  \draw[very thick](5.5, 0.2) -- (5,1);
\draw[fill=gray!30] (5,1) circle (1.5mm);
\draw[fill=gray!30] (6,1) circle (1.5mm);
\draw[fill=gray!30] (5.5, 0.2) circle (1.5mm);
 \draw[very thick](8,1) -- (9,1);
 \draw[very thick](8.5, 0.2) -- (8,1);
  \draw[very thick](8.5, 0.2) -- (9,1);
\draw[fill=gray!30] (8,1) circle (1.5mm);
\draw[fill=gray!30] (9,1) circle (1.5mm);
\draw[fill=gray!30] (8.5, 0.2) circle (1.5mm);
  \end{tikzpicture}
 \end{center}
\begin{center}
 \begin{tikzpicture}[scale=0.6]
     \node at (0,0.8) {\footnotesize \textit{(iii)}};
\draw[gray, ->] (-0.5,0.3) -- (0.5,0.3);
   \draw[thick](-1.25, 0.75) -- (-1.25,1.25) -- (-0.75,1.25) -- (-0.75,0.75) -- cycle;
   \draw[very thick] (-2,1) -- (-1,1);
 \draw[very thick](-1.5, 0.2) -- (-2,1);
  \draw[very thick](-1.5, 0.2) -- (-1,1);
\draw[fill=gray!30] (-2,1) circle (1.5mm);
\draw[fill=gray!30] (-1,1) circle (1.5mm);
\draw[fill=gray!30] (-1.5, 0.2) circle (1.5mm);
    \draw[very thick](1,1) -- (2,1);
  \draw[very thick](1.5, 0.2) -- (1,1);
    \draw[very thick](1.5, 0.2) -- (2,1);
    \draw[very thick](2,1) -- (2,3);
\draw[fill=red] (1,1) circle (1.5mm);
\draw[fill=red] (2,1) circle (1.5mm);
\draw[fill=red] (1.5, 0.2) circle (1.5mm);
    \draw[very thick](1,3) -- (2,3);
\draw[fill=gray!30] (1,3) circle (1.5mm);
\draw[fill=gray!30] (2,3) circle (1.5mm);
\draw[fill=gray!30] (1.5, 2.2) circle (1.5mm);
     \node at (7,0.8) {\footnotesize \textit{(iv)}};
\draw[gray, ->] (6.5,0.3) -- (7.5,0.3);
    \draw[thick](4.75, 0.75) -- (4.75,1.25) -- (5.25,1.25) -- (5.25,0.75) -- cycle;
    \draw[very thick](5,1) -- (6,1);
  \draw[very thick](5.5, 0.2) -- (5,1);
    \draw[very thick](5.5, 0.2) -- (6,1);
    \draw[very thick](6,1) -- (6,3);
\draw[fill=red] (5,1) circle (1.5mm);
\draw[fill=red] (6,1) circle (1.5mm);
\draw[fill=red] (5.5, 0.2) circle (1.5mm);
    \draw[very thick](5,3) -- (6,3);
\draw[fill=gray!30] (5,3) circle (1.5mm);
\draw[fill=gray!30] (6,3) circle (1.5mm);
\draw[fill=gray!30] (5.5, 2.2) circle (1.5mm);
    \draw[very thick](8,1) -- (9,1);
  \draw[very thick](8.5, 0.2) -- (8,1);
    \draw[very thick](8.5, 0.2) -- (9,1);
    \draw[very thick](9,1) -- (9,3);
    \draw[very thick](8,1) -- (8.5, 2.2);
\draw[fill=gray!30] (8,1) circle (1.5mm);
\draw[fill=red] (9,1) circle (1.5mm);
\draw[fill=red] (8.5, 0.2) circle (1.5mm);
    \draw[very thick](8,3) -- (9,3);
\draw[fill=gray!30] (8,3) circle (1.5mm);
\draw[fill=gray!30] (9,3) circle (1.5mm);
\draw[fill=gray!30] (8.5, 2.2) circle (1.5mm);
\node at (14,0.8) {\footnotesize \textit{(v)}};
\draw[gray, ->] (13.5,0.3) -- (14.5,0.3);
    \draw[thick](12.25, -0.05) -- (12.25,0.45) -- (12.75,0.45) -- (12.75,-0.05) -- cycle;
    \draw[very thick](12,1) -- (13,1);
  \draw[very thick](12.5, 0.2) -- (12,1);
    \draw[very thick](12.5, 0.2) -- (13,1);
    \draw[very thick](13,1) -- (13,3);
     \draw[very thick](12,1) -- (12.5, 2.2);
\draw[fill=gray!30] (12,1) circle (1.5mm);
\draw[fill=red] (13,1) circle (1.5mm);
\draw[fill=red] (12.5, 0.2) circle (1.5mm);
    \draw[very thick](12,3) -- (13,3);
\draw[fill=gray!30] (12,3) circle (1.5mm);
\draw[fill=gray!30] (13,3) circle (1.5mm);
\draw[fill=gray!30] (12.5, 2.2) circle (1.5mm);
    \draw[very thick](15,1) -- (16,1);
  \draw[very thick](15.5, 0.2) -- (15,1);
    \draw[very thick](15.5, 0.2) -- (16,1);
    \draw[very thick](16,1) -- (16,3);
    \draw[very thick](15,1) -- (15.5, 2.2);
    \draw[very thick](15.5,0.2) -- (15.5, 2.2);
\draw[fill=gray!30] (15,1) circle (1.5mm);
\draw[fill=red] (16,1) circle (1.5mm);
\draw[fill=gray!30] (15.5, 0.2) circle (1.5mm);
    \draw[very thick](15,3) -- (16,3);
\draw[fill=gray!30] (15,3) circle (1.5mm);
\draw[fill=gray!30] (16,3) circle (1.5mm);
\draw[fill=gray!30] (15.5, 2.2) circle (1.5mm);
\node at (21,0.8) {\footnotesize \textit{(vi)}};
\draw[gray, ->] (20.5,0.3) -- (21.5,0.3);
      \draw[thick](19.75, 0.75) -- (19.75,1.25) -- (20.25,1.25) -- (20.25,0.75) -- cycle;
    \draw[very thick](19,1) -- (20,1);
  \draw[very thick](19.5, 0.2) -- (19,1);
    \draw[very thick](19.5, 0.2) -- (20,1);
    \draw[very thick](20,1) -- (20,3);
     \draw[very thick](19,1) -- (19.5, 2.2);
      \draw[very thick](19.5,0.2) -- (19.5, 2.2);
\draw[fill=gray!30] (19,1) circle (1.5mm);
\draw[fill=red] (20,1) circle (1.5mm);
\draw[fill=gray!30] (19.5, 0.2) circle (1.5mm);
    \draw[very thick](19,3) -- (20,3);
\draw[fill=gray!30] (19,3) circle (1.5mm);
\draw[fill=gray!30] (20,3) circle (1.5mm);
\draw[fill=gray!30] (19.5, 2.2) circle (1.5mm);
\draw[very thick] (23,1) to [out=90, in=0] (22,3);
    \draw[very thick](22,1) -- (23,1);
  \draw[very thick](22.5, 0.2) -- (22,1);
    \draw[very thick](22.5, 0.2) -- (23,1);
    \draw[very thick](23,1) -- (23,3);
    \draw[very thick](22,1) -- (22.5, 2.2);
    \draw[very thick](22.5,0.2) -- (22.5, 2.2);
\draw[fill=gray!30] (22,1) circle (1.5mm);
\draw[fill=gray!30] (23,1) circle (1.5mm);
\draw[fill=gray!30] (22.5, 0.2) circle (1.5mm);
    \draw[very thick](22,3) -- (23,3);
\draw[fill=gray!30] (22,3) circle (1.5mm);
\draw[fill=gray!30] (23,3) circle (1.5mm);
\draw[fill=gray!30] (22.5, 2.2) circle (1.5mm);
  \end{tikzpicture}
 \end{center}
    \caption{(i) turns an $A^1$-structure into an $A^2$-structure, (ii) turns an $A^2$-structure into a triangle, (iii) assigns a triangle a responsible $A$-structure, (iv)-(vi) continue augmentation with a responsible $A$-structure. The squares indicate the offered vertices and the red color indicates which vertices in a triangle with responsibility still need to be offered one more time.}
    \label{fig:phase_2_simple}
\end{figure}

A triangle needs at least four hits (possibly more if they do not land on suitable vertices) until an augmentation happens, in which case the old triangle is removed and the two new triangles are created.
All other triangles that are in the middle of the augmentation process with respect to this augmented triangle are freed again and the edges between these triangles are ignored.
Note that it is also possible that an augmentable structure becomes a triangle before it is augmented, which is what we call a \emph{completion} and in which case also all edges to other triangles are ignored.
This phase takes roughly $2.259248n$ rounds and at its end we will have all but $10^{-15} n$ vertices in triangles.
Again we refer to \Cref{fig:plot} for an illustration of the development of the roles of vertices.
Observe that the abrupt change of the slope in the number of $A^2$-structures from Phase 1 to Phase 2 is due to the fact that an $A^1$-structure contains an unsaturated vertex whose offering creates an $A^2$-structure and hence this is more likely than hitting a matching edge.
Note that a strategy could use completion and augmentation from the beginning, but this would constitute great difficulties for its analysis and only for this reason, we resort to a separate first phase.

Similarly, the reason that the second phase stops before the triangle factor is complete is not because of the strategy, but because of the proof method.
We need a third clean-up phase, which is very similar to the second phase, but with a different analysis, and finishes the triangle factor in roughly $0.014739n$ rounds.

The attentive reader might have noted that the aforementioned durations sum up to $2.949498n$ rounds.
We obtain the improved bound of $2.929329$ in \Cref{main_theorem_k3_factor_upper} by already in Phase 1 using the hits on triangles to assign responsibilities (similarly to (iii) in \Cref{fig:phase_2_simple}) and prepare for later augmentation (see \Cref{fig:plot} to get an impression of the head start we get in Phase 2).

\subsection{Differential equation method}\label{sec:DEM}
Describing the strategy is only one side of the coin, but proving that it actually does something within a constrained time is the other.
For this we need to closely track the number of structures and types of augmentations in progress.
The perfect tool for this is the differential equation method (DEM) introduced by Wormald~\cite{wormald1995differential}.
In a nutshell what it does is when we want to track some variable, say the number of unsaturated vertices, we can turn the expected differences scaled by $n$ into a system of differential equations, whose solution in specific circumstances gives us w.h.p.~a good estimate of the truth.

More formally, for example, let $Y(t)$ be the number of vertices of $G$ that were not offered until step $t$ and $H(t)$ encode the history of the process.
Then $\mathbb{E}[Y(t+1)-Y(t) \mid H(t)]=-Y(t)/n$ and, scaling by $n$ and writing $x=t/n$, we transform this into the differential equation $y'(x)=-y(x)$ with initial condition $y(0)=1$ coming from $Y(0)=n$.
After verifying three hypotheses (boundedness, trend, and Lipschitz) the DEM-Theorem \cite[Theorem~5.1]{wormald1999differential} guarantees a unique solution such that $Y(xn)$ is w.h.p.~concentrated around $y(x)n=e^{-x}n$.
We will have a much more complex system, for which we sometimes can only find numerical solutions, but conceptually it remains the same idea.
We refer to a slightly stronger non-asymptotic version \cite[Theorem 2]{warnke2019wormald}, which is more suitable for our application and directly includes deviations from the initial conditions, and again to~\cite{eva_thesis} for a more in-depth discussion.

\subsection{Balls into Bins}

The sequence of vertices offered is closely related to the so-called \emph{Balls into Bins} model.
In this model $b$ unlabeled balls are thrown into $m$ labeled bins, where for every ball, one of the bins is chosen u.a.r.~and independently of all other choices.
Estimating how many vertices are offered $k$ times (see~\cite{ben2020semi}) is required to prove \Cref{main_theorem_k3_factor_offline}.
However, we will not rely on this model to track how many vertices are offered how often, but how many triangles are responsible for some augmentable structure as these are assigned u.a.r.\ as well.
In this case we need to control the maximum load of a bin, i.e.~the maximum number of triangles responsible for a given $A$-structure.

\begin{lm}\label{bib_max_load_lemma}
For $b = O(n)$ balls which are u.a.r.\ and independently assigned to  $m = \Theta(n)$ bins, the probability that a single bin gets at least $\log(n)$ balls is at most $O\left(\frac{1}{n^2}\right)$.
\end{lm}

The standard proof follows with Markov's inequality and we omit it here.
Due to augmentation or completion, augmentable structures will also disappear, which corresponds to the disappearance of a bin.
If the bin is selected u.a.r.~we expect $\tfrac{b}{m}$ balls to disappear alongside the bin.
To see this, note that after removing a bin the remaining balls are still distributed u.a.r.~and independently.
However, if we select a ball u.a.r.~and destroy the corresponding bin, large bins have a higher chance of being picked and, hence, the expected size of that bin is larger.
Actually, it is asymptotically one larger than $\tfrac{b}{m}$ as the next lemma shows.

\begin{lm}
\label{bib_destruction_lemma}
For $b = O(n)$ balls which are u.a.r.\ and independently assigned to  $m = \Theta(n)$ bins, suppose that one of the balls is chosen u.a.r.\ and its corresponding bin gets destroyed. Then the expected number of balls being destroyed is equal to
\begin{equation*}
    \frac{b-1}{m} + 1  =  \frac{b}{m} + 1 + O\left(\frac{1}{n}\right).
\end{equation*}
\end{lm}

This is also used implicitly in the proof for perfect matchings in \cite{gao2022perfect}, but  we feel it is worth stating it explicitly here.
For completeness we state the full proof which is also due to~\cite{pralat_personal_communication}.

\begin{proof}[Proof of \Cref{bib_destruction_lemma}]
We analyze the situation by deferring the decision about the placement of the balls until one of the $b$ balls is chosen.
This is sometimes called \emph{the principle of deferred decision}.
The ball whose bin will be destroyed is assigned u.a.r.~to one of the $m$ bins, say bin $j$.
The question is how many of the other $b-1$ balls are also in that bin.
In expectation this is $\frac{b-1}{m}$, so in total, for $m = \Theta(n)$, we expect bin $j$ to contain
\begin{equation*}
   1 + \frac{b-1}{m} = 1 + \frac{b}{m} - \frac{1}{m} =  1 + \frac{b}{m} + O\left(\frac{1}{n}\right)
\end{equation*}
balls.

For this to be feasible, it is crucial that only the question whether a ball goes to bin $j$ (and disappears along with it) was answered, and, otherwise, the assignment to the remaining $m-1$ bins remains u.a.r.~and can be further postponed.
This ensures that the placement of the remaining balls remains u.a.r.~and independent, which ensures that this procedure can be iterated while the underlying distribution is preserved.
\end{proof}

\subsection{Concentration inequalities}
\label{sec:concentration}
We additionally need the following standard concentration inequalities whose statements we take from~\cite{ben2020very}.
The first, a simplified version of Azuma's inequality, and the second, a simplified version of Talagrand's inequality.

\begin{lm}\label{lem:azuma}
Let $X$ be a non-negative random variable, not identically $0$, which is determined by $T$ independent trials. Assume that changing the outcome of any one of the trials can change the value of $X$ by at most $c >0$.

Then for every $\lambda \geq 0$
\begin{equation*}
    \mathbb{P} \left[ X \leq \mathbb{E}[X] - \lambda \right] \leq e^{- \frac{\lambda^2}{2\cdot c^2 \cdot T}} \hspace{0.4cm} \text{and} \hspace{0.4cm}  
    \mathbb{P} \left[ X \geq \mathbb{E}[X] + \lambda \right] \leq e^{- \frac{\lambda^2}{2\cdot c^2 \cdot T}}.
\end{equation*}
\end{lm}

\begin{lm}\label{lem:talagrand}
Let $X$ be a non-negative random variable, not identically $0$, which is determined by $T$ independent trials. Assume that changing the outcome of any one of the trials can change the value of $X$ by at most $c >0$.
Further assume that for every arbitrary $s > 0$ it holds that if $X\geq s$, then there is a set of at most $g\cdot s$ trials whose outcomes certify that $X\geq s$ (i.e. irrespective of the outcomes of all other trials, $X$ will remain at least $s$).
Then it holds for every $0 \leq \lambda \leq \mathbb{E}[X]$ that
\begin{equation*}
    \mathbb{P} \left[ |X - \mathbb{E}[X]| > \lambda + 60 \cdot c \cdot \sqrt{g \cdot \mathbb{E}[X]} \right] \leq 4 \cdot e^{- \frac{\lambda^2}{8\cdot c^2 \cdot g \cdot \mathbb{E}[X]}}.
\end{equation*}
\end{lm}

\section{Building a $K_3$-Factor in at most $2.929329 \, n$ rounds}
\label{k3_factor_section}

In this section we prove \Cref{main_theorem_k3_factor_upper}, that is, we will show that there exists a strategy $\mathcal{S}$, consisting of three phases, s.t.\ after $2.929329 n$ rounds w.h.p.~$G_t^{\mathcal{S}}$ contains a triangle factor.
This strategy works essentially as the one described in \Cref{strategy_description} with the additional ingredient that it already lays some groundwork for augmentation during the first phase by using the hits on complete triangles to assign responsibilities.

\subsection{Phase 1:  Creating the Disjoint Union of Augmentable and Complete Structures in at most $ 0.675511 \, n$ Rounds}\label{trianglefactor_phase1_section}

Recall the strategy for Phase~1 from \Cref{strategy_description}, where matching edges and $A^2$-structures are created.
The most natural idea how to include these responsibility assignments in Phase~1 would be to choose a matching edge or an $A^2$-structure u.a.r.~whenever a triangle without responsibility is offered.
However, given the varying times of ``birth'' of matching edges and $A^2$-structures in question and hence the varying number of rounds during which they would have had the opportunity to accumulate responsible triangles, an unpredictably uneven load (with different values for the average load of matching edges on the one hand and of $A^2$-structures on the other hand) would occur.

An obvious way to avoid that is to restrict this load to $1$ from the very beginning.
That is, whenever a vertex from a triangle without responsibility is offered, to assign the very same to a structure that does not have any responsible triangles yet.
This approach is most tractable if one explicitly distinguishes between $K_3$'s responsible for matching edges and $K_3$'s responsible for $A^2$-structures.
Our first goal is to arrive at a stage where $G_t$ contains the disjoint union of $A^1$-structures, $A^2$-structures, and complete triangles and, in addition, we would already like to assign responsibilities to some of the triangles.

 \begin{defi}\label{def_sets_phase_1_with_r3_k3}
At time $t \in \mathbb{N}$, let $\mathcal{M}_t$ be the set of matching edges, $\mathcal{A}^2_t$ be the set of $A^2$-structures, and $\mathcal{F}_t$ be the set of triangles.
Let $r_t : \mathcal{F}_t \rightarrow \{ \emptyset \} \cup \mathcal{M}_t \cup \mathcal{A}^2_t$ be the map that tracks the responsibilities, where $\emptyset$ represents ``no responsibility''. 
The vertex set $V(G_t)$ is partitioned into the following four sets:
\begin{itemize}[noitemsep]
    \item $F_t:=\cup_{F \in \mathcal{F}_t} V(F)$ --- the vertices in triangles
    \item $A^2_t:=\cup_{A \in \mathcal{A}^2_t} V(A)$ --- the vertices in $A^2$-structures
    \item $M_t:=\cup_{M \in \mathcal{M}_t} V(M)$ --- the vertices in matching edges
    \item $U_t:=V(G_t)\setminus (F_t \cup A^2_t \cup M_t)$ --- the untouched vertices.
\end{itemize}
Moreover, $F_t$ is partitioned into the following three sets:
\begin{itemize}[noitemsep]
\item $D_t:=\cup_{F \in r_t^{-1}(\emptyset)} V(F)$ --- the vertices in triangles without responsibility
\item $R^M_t:=\cup_{F \in r_t^{-1}(\mathcal{M}_t)} V(F)$ --- the vertices in triangles responsible for a matching edge
\item $R^A_t:=\cup_{F \in r_t^{-1}(\mathcal{A}^2_t)} V(F)$ --- the vertices in triangles responsible for an $A^2$-structure.
\end{itemize}
\end{defi}

Our algorithm will ensure that all graphs in $\mathcal{M}_t \cup \mathcal{A}^2_t \cup \mathcal{F}_t$ are subgraphs of $G_t$ and pairwise vertex-disjoint.
Recall that any responsibility $r_t(F)=A \not=\emptyset$ comes with an edge in $G_t$ from the triangle $F$ to a vertex in $A$, which can either be a matching edge or an $A^2$ structure.

\begin{figure}[htbp]
\begin{center}
 \begin{tikzpicture}[scale=0.8]
   \draw[thick](-8.25, 0.75) -- (-8.25,1.25) -- (-7.75,1.25) -- (-7.75,0.75) -- cycle;
\draw[fill=gray!30] (-8,1) circle (1.5mm);
\draw[fill=gray!30] (-8.5, 0.2) circle (1.5mm);
 \draw[very thick](-5.5, 0.2) -- (-5,1);
\draw[fill=gray!30] (-5,1) circle (1.5mm);
\draw[fill=gray!30] (-5.5, 0.2) circle (1.5mm);
   \draw[thick](-1.75, -0.05) -- (-1.75,0.45) -- (-1.25,0.45) -- (-1.25,-0.05) -- cycle;
 \draw[very thick](-1.5, 0.2) -- (-2,1);
\draw[fill=gray!30] (-2,1) circle (1.5mm);
\draw[fill=gray!30] (-1,1) circle (1.5mm);
\draw[fill=gray!30] (-1.5, 0.2) circle (1.5mm);
  \draw[very thick](1.5, 0.2) -- (1,1);
    \draw[very thick](1.5, 0.2) -- (2,1);
\draw[fill=gray!30] (1,1) circle (1.5mm);
\draw[fill=gray!30] (2,1) circle (1.5mm);
\draw[fill=gray!30] (1.5, 0.2) circle (1.5mm);
   \draw[thick,red](-2.4, 0.2) -- (-2.2, 0.5) -- (-2,0.2) -- cycle;
  \draw[thick, red](2, 0.2) -- (2.2, 0.5) -- (2.4,0.2) -- cycle;
    \draw[thick](4.75, 0.75) -- (4.75,1.25) -- (5.25,1.25) -- (5.25,0.75) -- cycle;
 \draw[very thick](5.5, 0.2) -- (6,1);
  \draw[very thick](5.5, 0.2) -- (5,1);
\draw[fill=gray!30] (5,1) circle (1.5mm);
\draw[fill=gray!30] (6,1) circle (1.5mm);
\draw[fill=gray!30] (5.5, 0.2) circle (1.5mm);
 \draw[very thick](8,1) -- (9,1);
 \draw[very thick](8.5, 0.2) -- (8,1);
  \draw[very thick](8.5, 0.2) -- (9,1);
\draw[fill=gray!30] (8,1) circle (1.5mm);
\draw[fill=gray!30] (9,1) circle (1.5mm);
\draw[fill=gray!30] (8.5, 0.2) circle (1.5mm);
   \draw[thick,red](4.6, 0.2) -- (4.8, 0.5) -- (5,0.2) -- cycle;
  \draw[thick](9, 0.2) -- (9.2, 0.5) -- (9.4,0.2) -- cycle;
   \node at (-6.8,0.8) {\footnotesize \ref{alg1_i}};
\draw[gray, ->] (-7.3,0.5) -- (-6.3,0.5);
   \node at (0,0.8) {\footnotesize {\ref{alg1_ii}}};
\draw[gray, ->] (-0.5,0.5) -- (0.5,0.5);
   \node at (7,0.8) {\footnotesize {\ref{alg1_iii}}};
\draw[gray, ->] (6.5,0.5) -- (7.5,0.5);
  \end{tikzpicture}
 \end{center}

\begin{center}
 \begin{tikzpicture}[scale=0.8]
     \node at (0,0.8) {\footnotesize {\ref{alg1_iv}}};
\draw[gray, ->] (-0.5,0.3) -- (0.5,0.3);
   \draw[thick](-1.25, 0.75) -- (-1.25,1.25) -- (-0.75,1.25) -- (-0.75,0.75) -- cycle;
   \draw[very thick] (-2,1) -- (-1,1);
 \draw[very thick](-1.5, 0.2) -- (-2,1);
  \draw[very thick](-1.5, 0.2) -- (-1,1);
\draw[fill=gray!30] (-2,1) circle (1.5mm);
\draw[fill=gray!30] (-1,1) circle (1.5mm);
\draw[fill=gray!30] (-1.5, 0.2) circle (1.5mm);
    \draw[very thick](1,1) -- (2,1);
  \draw[very thick](1.5, 0.2) -- (1,1);
    \draw[very thick](1.5, 0.2) -- (2,1);
    \draw[very thick](2,1) -- (2,2);
\draw[fill=red] (1,1) circle (1.5mm);
\draw[fill=red] (2,1) circle (1.5mm);
\draw[fill=red] (1.5, 0.2) circle (1.5mm);
    \draw[very thick](1,2) -- (2,2);
\draw[fill=gray!30] (1,2) circle (1.5mm);
\draw[fill=gray!30] (2,2) circle (1.5mm);
  \end{tikzpicture}
 \end{center}
    \caption{The four relevant possibilities in Phase~1. The squares indicate the offered vertices $v_{t+1}$.
    The small triangles indicate responsible triangles for that structure (red) or no longer responsible ones (black).
    Each red vertex still needs to be offered one more time for an augmentation, but that will only be tracked in Phase~2.}
    \label{fig:phase_1}
\end{figure}
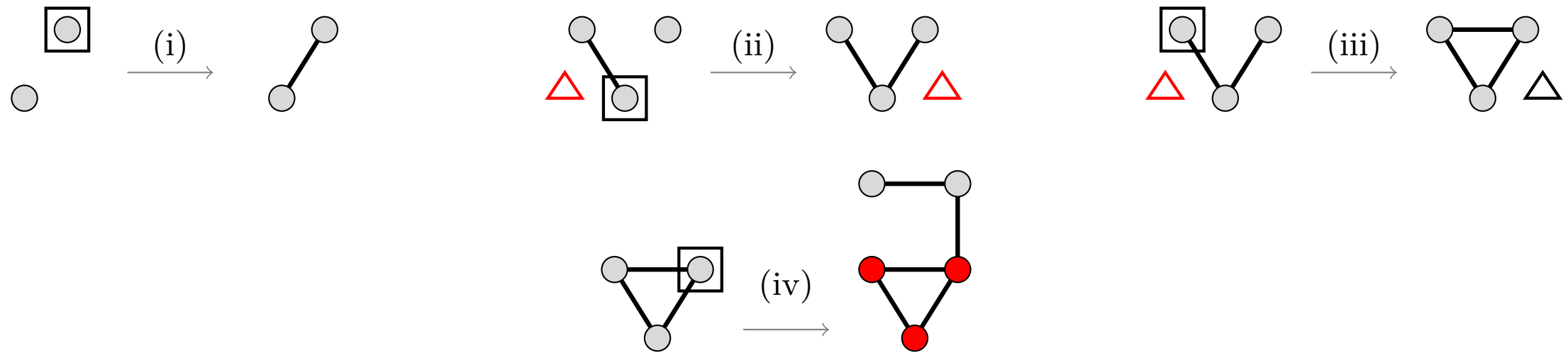

Now we provide the formal details of our algorithm for Phase~1 that is also illustrated in \Cref{fig:phase_1}.
In each update step we will only mention the changes for the sets, i.e.~$\mathcal{M}_{t+1}:=\mathcal{M}_t$, $\mathcal{A}^2_{t+1}:=\mathcal{A}^2_t$, $\mathcal{F}_{t+1}:=\mathcal{F}_t$, and $r_{t+1}(F):=r_t(F)$ for any $F \in \mathcal{F}_{t+1}$ if not specified otherwise.
The update of the other sets is automatic according to \Cref{def_sets_phase_1_with_r3_k3}.

 \begin{alg}\label{strat_phase_1_with_r3_k3}
 Initialize $\mathcal{M}_0=\mathcal{F}_0=\mathcal{A}^2_0=\emptyset$ (this gives $U_0 = [n]$ and $M_0 = A^2_0 = F_0 = D_0 = R^M_0 = R^A_0 = \emptyset$).
 In round $t \ge 0$ if $|U_t| \le |\mathcal{M}_t|$ terminate and, otherwise, let $v_{t+1}$ be the offered vertex.
\begin{enumerate}[label=(\roman*), itemsep=0.7ex]
\item \label{alg1_i} If $v_{t+1} \in U_{t}$, choose $w_{t+1}$ u.a.r.\ from $ U_{t} \setminus \{v_{t+1}\}$.
Let $\mathcal{M}_{t+1}:=\mathcal{M}_t \cup \{ \{v_{t+1}w_{t+1} \} \}$.
If $U_{t} \setminus \{v_{t+1}\}=\emptyset$ abort the algorithm.
\item \label{alg1_ii} If $v_{t+1} \in M$ for some $M=\{v_{t+1}u\} \in \mathcal{M}_t$, choose $w_{t+1}$ u.a.r.\ from $ U_{t}$ and create a new $A^2$-structure $A=\{ w_{t+1}v_{t+1}, v_{t+1}u\}$.
We set $\mathcal{M}_{t+1}:=\mathcal{M}_t \setminus \{M\}$, $\mathcal{A}^2_{t+1}:=\mathcal{A}^2_t \cup \{ A\}$, and if there exists $H \in r_t^{-1}(M)$, then $r_{t+1}(H):=A$.
If $U_{t}=\emptyset$ abort the algorithm.
\item \label{alg1_iii} If $v_{t+1}$ is an endpoint of some $A^2$-structure $A=\{v_{t+1}u_1,u_1u_2\} \in \mathcal{A}^2_t$, we set $w_{t+1}:=u_2$ and thus create a new triangle $F = A\cup \{ v_{t+1}w_{t+1} \}$.
We set $\mathcal{A}^2_{t+1}:=\mathcal{A}^2_t \setminus \{A\}$, $\mathcal{F}_{t+1}:=\mathcal{F}_t \cup \{ F\}$, $r_{t+1}(F):=\emptyset$, and if there exists  $H \in r_t^{-1}(A)$, then $r_{t+1}(H):=\emptyset$. 
\item \label{alg1_iv} If $v_{t+1} \in F$ for some $F \in \mathcal{F}_t$ with $r_t(F) = \emptyset$, choose a matching edge $M= \{w_{t+1}u\}$ u.a.r.~from  $\mathcal{M}_t \setminus r_t(\mathcal{F}_t)$ and set $r_{t+1}(F):=M$.
If $\mathcal{M}_t \setminus r_t(\mathcal{F}_t) = \emptyset$ abort the algorithm.
\item \label{alg1_v} Else make an arbitrary move and forget about it.
\end{enumerate}
\end{alg}

 Recall that we assume $n$ to be divisible by three and, hence, the algorithm terminates as soon as $|U_t| = |\mathcal{M}_t|$.
 Note that an $A^1$-structure consists of exactly one matching edge (whose number is equal to $|\mathcal{M}_t|$) and exactly one isolated vertex (whose number is equal to $|U_t|$), so at the point of termination $G_t$ does contain the disjoint union of augmentable and complete structures.
 We claim that w.h.p.~the algorithm terminates with $|U_t| = |\mathcal{M}_t|$.
 
\begin{prp}
\label{running_time_phase_1_with_r3_k3}
    \Cref{strat_phase_1_with_r3_k3} w.h.p.\ terminates after $t_1 \le 0.675511 n$ rounds with $|U_{t_1}| = |\mathcal{M}_{t_1}|$.
    Moreover, setting $s \coloneqq \frac{\sqrt{30}}{4}\cdot \arctan(\sqrt{2})$ it holds w.h.p.\ that
\begin{align*}
|U_{t_1}| &= \frac{1}{\sqrt{3}\cdot e^{\sqrt{2} \cdot \arctan(\sqrt{2})}} \cdot n + O(n^{\frac{3}{4}})\\
|M_{t_1}|  &= \frac{2}{\sqrt{3}\cdot e^{\sqrt{2} \cdot \arctan(\sqrt{2})}} \cdot n + O(n^{\frac{3}{4}})\\
|A^2_{t_1}| &= \frac{3 - \sqrt{3}}{ e^{\sqrt{2} \cdot \arctan(\sqrt{2})}} \cdot n +  O(n^{\frac{3}{4}})\\
|F_{t_1}| &=  \left(1 - \frac{3}{e^{\sqrt{2} \cdot \arctan(\sqrt{2})}} \right) \cdot n+  O(n^{\frac{3}{4}})\\
|D_{t_1}| &= \left( \frac{1}{4} - \frac{\frac{18}{17}}{\sqrt{3}\cdot e^{\sqrt{2} \cdot \arctan(\sqrt{2})}} 
+  \frac{\frac{39}{68} \cdot \cos( s) + \frac{27\sqrt{15}}{340} \cdot \sin( s)}{e^{\frac{7\sqrt{2}}{4} \cdot \arctan(\sqrt{2})}} \right)
\cdot n + O(n^{\frac{3}{4}})\\
|R^M_{t_1}| &= \left( \frac{3}{8} - \frac{\frac{36}{17}}{\sqrt{3}\cdot e^{\sqrt{2} \cdot \arctan(\sqrt{2})}} 
  - \frac{\frac{99}{136} \cdot \cos(s) - \frac{57\sqrt{15}}{680} \cdot \sin( s)}{e^{\frac{7\sqrt{2}}{4} \cdot \arctan(\sqrt{2})}} \right)
\cdot n + O(n^{\frac{3}{4}})\\
|R^A_{t_1}| &= \left( \frac{3}{8} + \frac{\frac{54}{17} - 3\sqrt{3}}{\sqrt{3}\cdot e^{\sqrt{2} \cdot \arctan(\sqrt{2})}} 
+ \frac{\frac{21}{136} \cdot \cos(s) - \frac{111\sqrt{15}}{680} \cdot \sin( s)}{e^{\frac{7\sqrt{2}}{4} \cdot \arctan(\sqrt{2})}} \right)
\cdot n + O(n^{\frac{3}{4}}).
\end{align*}
\end{prp}

See \Cref{fig:plot} for a plot of the functions throughout Phase 1.
Note that it directly follows from the definition of the algorithm that indeed the elements in $\mathcal{M}_{t_1}$, $\mathcal{A}^2_{t_1}$, and $\mathcal{F}_{t_1}$ are all pairwise vertex-disjoint subgraphs of $G_{t_1}$ and that for $F \in \mathcal{F}_{t_1}$ with $r_{t_1}(F) =A\not= \emptyset$ the edge $f_3a_1$  for some $f_3 \in V(F)$ and $a_1a_2 \in E(A)$ is also in $G_{t_1}$.

\begin{proof}
We will prove the proposition using the DEM, as explained in \Cref{sec:DEM}.
For $t\ge0$ we define the following random variables analogous to the vertex sets defined in \Cref{def_sets_phase_1_with_r3_k3}:
$U(t) :=|U_t|$, $M(t):=|M_t|$, $A^2(t):=|A^2_t|$, $F(t):=|F_t|$, $D(t):=|D_t|$, $R^M(t):=|R^M_t|$, and $R^A(t):=|R^A_t|$.
Following the initialization of \Cref{running_time_phase_1_with_r3_k3} it holds that $U(0) = n$ and $M(0) = A^2(0) = F(0) = D(0) = R^M(0)=  R^A(0)= 0$.
We need to understand the expected differences, which model the expected change from time step $t$ to $t+1$, for each of the random variables.
The variable $U(t)$ only changes in rule \ref{alg1_i}, where it decreases by two, and in rule \ref{alg1_ii}, where it decreases by one.
Hence,
\begin{equation}
\label{eq: Phase I U_t}
\begin{aligned}
\mathbb{E}[U(t+1) - U(t) \mid H(t)] =&  \mathbb{P}[v_{t+1} \in U(t) \mid H(t)] \cdot (-2) \\
&+ \mathbb{P}[v_{t+1} \in M(t) \mid H(t)] \cdot (-1)\\
=& \dfrac{U(t)}{n} 
\cdot (-2)  +  \dfrac{M(t)}{n}\cdot (-1) \\
\end{aligned}
\end{equation}

Similarly, we derive the next three directly from the definition of the algorithm and rules \ref{alg1_i}, \ref{alg1_ii}, and \ref{alg1_iii}
\begin{equation}
\label{eq: Phase I M_t}
\begin{aligned}
  \mathbb{E}[M(t+1) - M(t) \mid H(t)]&=
  \dfrac{U(t)}{n} \cdot 2 +  \dfrac{M(t)}{n}\cdot (-2),\\
  \mathbb{E}[A^2(t+1) - A^2(t) \mid H(t)]&=
\dfrac{M(t)}{n}  \cdot 3  +   \dfrac{\frac{2}{3} \cdot A^2(t)}{n} \cdot (-3),\\
  \mathbb{E}[F(t+1) - F(t) \mid H(t)] &= 
  \dfrac{\frac{2}{3} \cdot A^2(t)}{n} \cdot 3,
\end{aligned}
\end{equation}

For the next three we first need to understand some more complex probabilities related to responsible triangles
\begin{align*}
 \mathbb{P}\left[ v_{t+1} \in A^2_t \wedge \{v_{t+1}u_1,u_1u_2\} \in r(\mathcal{F}_t)\right]  &= 
 \dfrac{\frac{2}{3} \cdot R^A(t)}{n},\\
\mathbb{P}\left[ v_{t+1} \in A^2_t \wedge \{v_{t+1}u_1,u_1u_2\} \not\in r(\mathcal{F}_t)\right]  &= 
\dfrac{\frac{2}{3} \cdot (A^2(t) -  R^A(t))}{n},\\
\mathbb{P}\left[ v_{t+1} \in M_{t} \wedge \{v_{t+1}u\} \in r(\mathcal{F}_t)\right] &=
 \dfrac{\frac{2}{3} \cdot R^M(t)}{n},
\end{align*}

and then get
\begin{equation}
\label{eq: Phase I D_t}
\begin{aligned}
    \mathbb{E}[D(t+1) - D(t) \mid H(t)] &= \dfrac{D(t)}{n} \cdot (-3) 
+ \dfrac{\frac{2}{3} \cdot R^A(t)}{n}  \cdot 6
+ \dfrac{\frac{2}{3} \cdot (A^2(t) -  R^A(t))}{n} \cdot 3,\\
    \mathbb{E}[R^M(t+1) - R^M(t) \mid H(t)] &= 
    \dfrac{D(t)}{n} \cdot 3 
     +  \dfrac{\frac{2}{3} \cdot R^M(t)}{n} \cdot (-3) , \\
    \mathbb{E}[R^A(t+1) - R^A(t) \mid H(t)] &= 
  \dfrac{\frac{2}{3} \cdot R^M(t)}{n}  \cdot 3 
    + \dfrac{\frac{2}{3} \cdot R^A(t)}{n} \cdot (-3).
\end{aligned}
\end{equation}

From the expected differences in \Cref{eq: Phase I D_t,eq: Phase I M_t,eq: Phase I U_t} we create a system of differential equations by scaling $x = \frac{t}{n}$.

\begin{minipage}{0.5\textwidth}
\begin{center}
       \begin{align*}
   u'(x) &= -2\cdot u(x) - m(x) \\
   m'(x) &= 2\cdot u(x) - 2 \cdot m(x)\\
   a^{2\prime}(x) &= 3\cdot m(x) -2 \cdot a^2(x)\\
   f'(x) &= 2\cdot a^2(x)\\
   d'(x) &= -3\cdot d(x) + 2 \cdot r^A(x) + 2\cdot a^2(x) \\
   r^{M\prime}(x) &= 3\cdot d(x) - 2 \cdot r^M(x)\\
   r^{A\prime}(x) &= 2 \cdot r^M(x) - 2 \cdot r^A(x)
\end{align*}

\end{center}
\end{minipage}\hspace{2em}
\begin{minipage}{0.2\textwidth}
      \begin{align*}
   u(0) &= 1\\
   m(0) &= 0\\
   a^2(0) &= 0\\
   f(0) &= 0\\
   d(0) &= 0\\
   r^M(0) &= 0\\
   r^A(0) &= 0
\end{align*}
\end{minipage}

\bigskip

As the boundedness and trend hypotheses ($\beta(n) = 6$, $\gamma(n) = 0$, $\lambda_1(n) = 0$), as well as the Lipschitz hypothesis (the right hand sides constitute linear functions) are fulfilled, we can conclude with \cite[Theorem~5.1]{wormald1999differential} that this system has a unique solution.
The solution is
\begin{equation}
\label{eq: Phase I solution}
\begin{aligned}
   u(x) &= \frac{\cos(\sqrt{2}x)}{e^{2x}}\\
   m(x) &= \frac{ \sqrt{2} \cdot \sin(\sqrt{2}x)}{e^{2x}}\\
   a^2(x) &= \frac{ 3- 3 \cdot \cos(\sqrt{2}x)}{e^{2x}}\\
   f(x) &= \frac{ 2 \cdot \cos(\sqrt{2}x)- \sqrt{2} \cdot \sin(\sqrt{2}x) - 3}{e^{2x}} + 1\\
   d(x) &= \frac{1}{4} 
   - \frac{\frac{14}{17} \cdot \cos(\sqrt{2}x) + \frac{2\sqrt{2}}{17} \cdot \sin(\sqrt{2}x)}{e^{2x}}
+  \frac{\frac{39}{68} \cdot \cos( \frac{\sqrt{15}}{2}x) + \frac{27\sqrt{15}}{340} \cdot \sin( \frac{\sqrt{15}}{2}x)}{e^{3.5x}} \\
   r^M(x) &= \frac{3}{8} 
   + \frac{\frac{6}{17} \cdot \cos(\sqrt{2}x) - \frac{21\sqrt{2}}{17} \cdot \sin(\sqrt{2}x)}{e^{2x}}
   - \frac{\frac{99}{136} \cdot \cos( \frac{\sqrt{15}}{2}x) - \frac{57\sqrt{15}}{680} \cdot \sin( \frac{\sqrt{15}}{2}x)}{e^{3.5x}}\\
   r^A(x) &= \frac{3}{8} 
   + \frac{\frac{42}{17} \cdot \cos(\sqrt{2}x) + \frac{6\sqrt{2}}{17} \cdot \sin(\sqrt{2}x) - 3}{e^{2x}}
   + \frac{\frac{21}{136} \cdot \cos( \frac{\sqrt{15}}{2}x) - \frac{111\sqrt{15}}{680} \cdot \sin( \frac{\sqrt{15}}{2}x)}{e^{3.5x}}.
\end{aligned}
\end{equation}

Choosing $\lambda(n) \coloneqq n^{-\frac{1}{4}}$ yields that $U(x n)$ is w.h.p.\ concentrated around $e^{-2x}\cos(\sqrt{2}x) n$ with a deviation of at most $O\left(n^{\frac{3}{4}}\right)$ and analogously for the remaining random variables.
To estimate if and when \Cref{strat_phase_1_with_r3_k3} terminates, we look for the smallest $x$ such that $u(x) \leq \frac{1}{2} \cdot m(x)$, which implies
\begin{align*}
 \frac{\cos(\sqrt{2}x)}{e^{2x}} &\leq \frac{1}{2} \cdot \frac{ \sqrt{2} \cdot \sin(\sqrt{2}x)}{e^{2x}}\\
x &\geq \frac{\arctan(\sqrt{2}) }{\sqrt{2}} = 0.6755108588560399...
\end{align*}

Hence, after $t=0.675511 n$ rounds, $U(t) \leq \frac{1}{2} \cdot M(t)$ w.h.p.\ and \Cref{strat_phase_1_with_r3_k3} has terminated. 
Moreover, the analytical solution in \Cref{eq: Phase I solution} shows that w.h.p.~for all $x \in \left[0,\arctan(\sqrt{2})/\sqrt{2}\right]$ we have $r^M(x) \leq m(x)$ and when $r^M(x) > 0$ this implies $\frac{r^M(x)}{3} < \frac{m(x)}{2}$ and $\mathcal{M}_t \setminus r_t(\mathcal{F}_t) \not= \emptyset$.
Therefore, \Cref{strat_phase_1_with_r3_k3} does not abort because of rule \ref{alg1_iv}.
Similarly, $u(x)>0$ in this regime and the algorithm also does not abort in rules \ref{alg1_i} or \ref{alg1_ii} because throughout we have $|U_{t}| > 1$.
\end{proof}

\subsection{Phase II: Arriving at $(1-\alpha)\, \frac{n}{3}$ $K_3$'s in $ 2.251183 \, n$ more rounds}\label{subsection_k3_phase2}

So far, the SRGP went on for $t_1 \leq 0.675511\cdot n$ rounds and $\left(1 - \frac{3}{e^{\sqrt{2} \cdot \arctan(\sqrt{2})}} \right) \cdot n \approx 0.223074 \cdot n$ of the vertices of $G_t$ are w.h.p.\ contained in complete triangles, while all remaining vertices are contained in augmentable structures.
The goal is now to ensure that for $\alpha \coloneqq \frac{1}{10^{21}}$ at the end of Phase II, $(1-\alpha)\cdot n$ vertices will be contained in triangles. To this end, both  completion and augmentation are used.
We emphasize that despite the more complex definitions and larger number of variables we proceed analogously to Phase I as conceptually it is the same argument.
After defining all relevant sets and the algorithm, we prove its termination with the DEM.

Like in Phase 1, we have $\mathcal{A}^2_t$ and $\mathcal{F}_t$ and additionally a set of $A^1$-structures $\mathcal{A}^1_t$.
At the end of Phase 1 we have $|\mathcal{M}_t|=|U_t|$ so we can order both sets and initialize  $\mathcal{A}^1_0$ by combining a matching edge $a_1a_2$ with an untouched vertex $a_3$ into the $A^1$-structure $A$ with $V(A) = \{a_1, a_2,a_3\}$ and $E(A) = \{a_1a_2\}$. 
We will still have responsibilities of the triangles $\mathcal{F}_t$ encoded by $r_t$.
Moreover, these triangles go through different stages of augmentation where 3, 2, 1 vertices still need to be hit (cf.~\Cref{fig:augmentation_triangles}) and we will encode this with another map $\rho_t$.

One additional complication that arises from Phase I is that there is a higher proportion of $A^1$-structures which have a responsible triangle ($\approx 16 \%$ of them) than there is of $A^2$-structures ($\approx 6 \%$ of them).
To account for this, we conduct the analysis of Phase 2 in such a way that an explicit distinction between responsible triangles for $A^1$- and $A^2$-structures is maintained. This necessarily leads to separate random variables tracking the development of the respective responsible triangles.

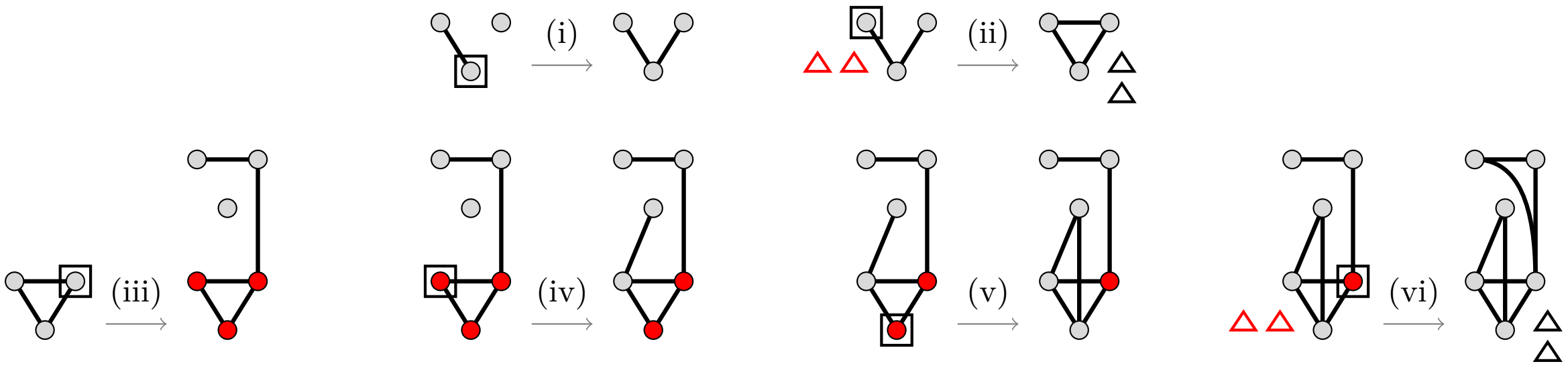
\begin{figure}[htbp]
 \begin{center}
 \begin{tikzpicture}[scale=0.6]
     \node at (0,0.8) {\footnotesize \ref{alg2_i}};
\draw[gray, ->] (-0.5,0.3) -- (0.5,0.3);
   \draw[thick](-1.75, -0.05) -- (-1.75,0.45) -- (-1.25,0.45) -- (-1.25,-0.05) -- cycle;
 \draw[very thick](-1.5, 0.2) -- (-2,1);
\draw[fill=gray!30] (-2,1) circle (1.5mm);
\draw[fill=gray!30] (-1,1) circle (1.5mm);
\draw[fill=gray!30] (-1.5, 0.2) circle (1.5mm);
  \draw[very thick](1.5, 0.2) -- (1,1);
    \draw[very thick](1.5, 0.2) -- (2,1);
\draw[fill=gray!30] (1,1) circle (1.5mm);
\draw[fill=gray!30] (2,1) circle (1.5mm);
\draw[fill=gray!30] (1.5, 0.2) circle (1.5mm);
     \node at (7,0.8) {\footnotesize \ref{alg2_ii}};
\draw[gray, ->] (6.5,0.3) -- (7.5,0.3);
    \draw[thick](4.75, 0.75) -- (4.75,1.25) -- (5.25,1.25) -- (5.25,0.75) -- cycle;
 \draw[very thick](5.5, 0.2) -- (6,1);
  \draw[very thick](5.5, 0.2) -- (5,1);
\draw[fill=gray!30] (5,1) circle (1.5mm);
\draw[fill=gray!30] (6,1) circle (1.5mm);
\draw[fill=gray!30] (5.5, 0.2) circle (1.5mm);
  \draw[thick,red](4, 0.2) -- (4.2, 0.5) -- (4.4,0.2) -- cycle;
   \draw[thick,red](4.6, 0.2) -- (4.8, 0.5) -- (5,0.2) -- cycle;
  \draw[thick](9, 0.2) -- (9.2, 0.5) -- (9.4,0.2) -- cycle;
\draw[thick](9, -0.3) -- (9.2, 0) -- (9.4, -0.3) -- cycle;
 \draw[very thick](8,1) -- (9,1);
 \draw[very thick](8.5, 0.2) -- (8,1);
  \draw[very thick](8.5, 0.2) -- (9,1);
\draw[fill=gray!30] (8,1) circle (1.5mm);
\draw[fill=gray!30] (9,1) circle (1.5mm);
\draw[fill=gray!30] (8.5, 0.2) circle (1.5mm);
  \end{tikzpicture}
 \end{center}
 
\begin{center}
 \begin{tikzpicture}[scale=0.6]
     \node at (0,0.8) {\footnotesize \ref{alg2_iii}};
\draw[gray, ->] (-0.5,0.3) -- (0.5,0.3);
   \draw[thick](-1.25, 0.75) -- (-1.25,1.25) -- (-0.75,1.25) -- (-0.75,0.75) -- cycle;
   \draw[very thick] (-2,1) -- (-1,1);
 \draw[very thick](-1.5, 0.2) -- (-2,1);
  \draw[very thick](-1.5, 0.2) -- (-1,1);
\draw[fill=gray!30] (-2,1) circle (1.5mm);
\draw[fill=gray!30] (-1,1) circle (1.5mm);
\draw[fill=gray!30] (-1.5, 0.2) circle (1.5mm);
    \draw[very thick](1,1) -- (2,1);
  \draw[very thick](1.5, 0.2) -- (1,1);
    \draw[very thick](1.5, 0.2) -- (2,1);
    \draw[very thick](2,1) -- (2,3);
\draw[fill=red] (1,1) circle (1.5mm);
\draw[fill=red] (2,1) circle (1.5mm);
\draw[fill=red] (1.5, 0.2) circle (1.5mm);
    \draw[very thick](1,3) -- (2,3);
\draw[fill=gray!30] (1,3) circle (1.5mm);
\draw[fill=gray!30] (2,3) circle (1.5mm);
\draw[fill=gray!30] (1.5, 2.2) circle (1.5mm);
     \node at (7,0.8) {\footnotesize \ref{alg2_iv}};
\draw[gray, ->] (6.5,0.3) -- (7.5,0.3);
    \draw[thick](4.75, 0.75) -- (4.75,1.25) -- (5.25,1.25) -- (5.25,0.75) -- cycle;
    \draw[very thick](5,1) -- (6,1);
  \draw[very thick](5.5, 0.2) -- (5,1);
    \draw[very thick](5.5, 0.2) -- (6,1);
    \draw[very thick](6,1) -- (6,3);
\draw[fill=red] (5,1) circle (1.5mm);
\draw[fill=red] (6,1) circle (1.5mm);
\draw[fill=red] (5.5, 0.2) circle (1.5mm);
    \draw[very thick](5,3) -- (6,3);
\draw[fill=gray!30] (5,3) circle (1.5mm);
\draw[fill=gray!30] (6,3) circle (1.5mm);
\draw[fill=gray!30] (5.5, 2.2) circle (1.5mm);
    \draw[very thick](8,1) -- (9,1);
  \draw[very thick](8.5, 0.2) -- (8,1);
    \draw[very thick](8.5, 0.2) -- (9,1);
    \draw[very thick](9,1) -- (9,3);
    \draw[very thick](8,1) -- (8.5, 2.2);
\draw[fill=gray!30] (8,1) circle (1.5mm);
\draw[fill=red] (9,1) circle (1.5mm);
\draw[fill=red] (8.5, 0.2) circle (1.5mm);
    \draw[very thick](8,3) -- (9,3);
\draw[fill=gray!30] (8,3) circle (1.5mm);
\draw[fill=gray!30] (9,3) circle (1.5mm);
\draw[fill=gray!30] (8.5, 2.2) circle (1.5mm);
     \node at (14,0.8) {\footnotesize \ref{alg2_v}};
\draw[gray, ->] (13.5,0.3) -- (14.5,0.3);
    \draw[thick](12.25, -0.05) -- (12.25,0.45) -- (12.75,0.45) -- (12.75,-0.05) -- cycle;
    \draw[very thick](12,1) -- (13,1);
  \draw[very thick](12.5, 0.2) -- (12,1);
    \draw[very thick](12.5, 0.2) -- (13,1);
    \draw[very thick](13,1) -- (13,3);
     \draw[very thick](12,1) -- (12.5, 2.2);
\draw[fill=gray!30] (12,1) circle (1.5mm);
\draw[fill=red] (13,1) circle (1.5mm);
\draw[fill=red] (12.5, 0.2) circle (1.5mm);
    \draw[very thick](12,3) -- (13,3);
\draw[fill=gray!30] (12,3) circle (1.5mm);
\draw[fill=gray!30] (13,3) circle (1.5mm);
\draw[fill=gray!30] (12.5, 2.2) circle (1.5mm);
    \draw[very thick](15,1) -- (16,1);
  \draw[very thick](15.5, 0.2) -- (15,1);
    \draw[very thick](15.5, 0.2) -- (16,1);
    \draw[very thick](16,1) -- (16,3);
    \draw[very thick](15,1) -- (15.5, 2.2);
    \draw[very thick](15.5,0.2) -- (15.5, 2.2);
\draw[fill=gray!30] (15,1) circle (1.5mm);
\draw[fill=red] (16,1) circle (1.5mm);
\draw[fill=gray!30] (15.5, 0.2) circle (1.5mm);
    \draw[very thick](15,3) -- (16,3);
\draw[fill=gray!30] (15,3) circle (1.5mm);
\draw[fill=gray!30] (16,3) circle (1.5mm);
\draw[fill=gray!30] (15.5, 2.2) circle (1.5mm);
     \node at (21,0.8) {\footnotesize \ref{alg2_vi}};
\draw[gray, ->] (20.5,0.3) -- (21.5,0.3);
      \draw[thick](19.75, 0.75) -- (19.75,1.25) -- (20.25,1.25) -- (20.25,0.75) -- cycle;
    \draw[very thick](19,1) -- (20,1);
  \draw[very thick](19.5, 0.2) -- (19,1);
    \draw[very thick](19.5, 0.2) -- (20,1);
    \draw[very thick](20,1) -- (20,3);
     \draw[very thick](19,1) -- (19.5, 2.2);
      \draw[very thick](19.5,0.2) -- (19.5, 2.2);
\draw[fill=gray!30] (19,1) circle (1.5mm);
\draw[fill=red] (20,1) circle (1.5mm);
\draw[fill=gray!30] (19.5, 0.2) circle (1.5mm);
    \draw[very thick](19,3) -- (20,3);
\draw[fill=gray!30] (19,3) circle (1.5mm);
\draw[fill=gray!30] (20,3) circle (1.5mm);
\draw[fill=gray!30] (19.5, 2.2) circle (1.5mm);
\draw[very thick] (23,1) to [out=90, in=0] (22,3);
    \draw[very thick](22,1) -- (23,1);
  \draw[very thick](22.5, 0.2) -- (22,1);
    \draw[very thick](22.5, 0.2) -- (23,1);
    \draw[very thick](23,1) -- (23,3);
    \draw[very thick](22,1) -- (22.5, 2.2);
    \draw[very thick](22.5,0.2) -- (22.5, 2.2);
\draw[fill=gray!30] (22,1) circle (1.5mm);
\draw[fill=gray!30] (23,1) circle (1.5mm);
\draw[fill=gray!30] (22.5, 0.2) circle (1.5mm);
    \draw[very thick](22,3) -- (23,3);
\draw[fill=gray!30] (22,3) circle (1.5mm);
\draw[fill=gray!30] (23,3) circle (1.5mm);
\draw[fill=gray!30] (22.5, 2.2) circle (1.5mm);
  \draw[thick,red](18, 0.2) -- (18.2, 0.5) -- (18.4,0.2) -- cycle;
   \draw[thick,red](18.6, 0.2) -- (18.8, 0.5) -- (19,0.2) -- cycle;
  \draw[thick](23, 0.2) -- (23.2, 0.5) -- (23.4,0.2) -- cycle;
   \draw[thick](23, -0.3) -- (23.2, 0) -- (23.4,-0.3) -- cycle;
   \end{tikzpicture}
 \end{center}
    \caption{
    Illustration of the steps of \Cref{strat_phase_2_with_distinction_k3}.
    The squares indicate the offered vertices $v_{t+1}$.
    For a triangle $F$ the vertices $\rho_t(F)$, which still need one more hit, are colored red.
    The small triangles indicate responsible triangles for that structure (red) or no longer responsible ones (black).
    Note that rules \ref{alg2_iii}-\ref{alg2_vi} depict responsibilities for an $A^1$-structure, but the procedure is identical for an $A^2$-structure. Also in rules \ref{alg2_iv}-\ref{alg2_vi} the picture might look different depending on when the vertex offered in rule \ref{alg2_iii} is offered again.}
    \label{fig:phase_2_again}
\end{figure}

 \begin{defi}\label{def_sets_phase_2_with_distinction_k3}
At time $t \in \mathbb{N}$, let $\mathcal{A}^1_t$ be the set of $A^1$-structures, $\mathcal{A}^2_t$ be the set of $A^2$-structures, and $\mathcal{F}_t$ be the set of triangles.
Let $r_t : \mathcal{F}_t \rightarrow \{ \emptyset \} \cup \mathcal{A}^1_t \cup \mathcal{A}^2_t$ be the map that tracks the responsibilities, where $\emptyset$ represents ``no responsibility''. 
The vertex set $V(G_t)$ is partitioned into the following three sets:
\begin{itemize}[noitemsep]
    \item $F_t:=\cup_{F \in \mathcal{F}_t} V(F)$ --- the vertices in triangles
    \item $A^2_t:=\cup_{A \in \mathcal{A}^2_t} V(A)$ --- the vertices in $A^2$-structures
    \item $A^1_t:=\cup_{A \in \mathcal{A}^1_t} V(A)$ --- the vertices in $A^1$-structures.
\end{itemize}
Let $\rho_t : \mathcal{F}_t \rightarrow 2^{F_t}$ be the map that identifies those vertices of a triangle that still need to be offered in order for an augmentation to be successful, where $\emptyset$ represents ``no responsibility''.
\end{defi}
Note that in accordance with the augmentation procedure we will ensure that $|\rho_t(F)| \in \{0,1,2,3\}$ and that $r_t(F) = \emptyset \Leftrightarrow \rho_t(F) = \emptyset$.
Now we provide the formal details of our algorithm for Phase~2 that is illustrated in~\Cref{fig:phase_2_again}.
In each update step, it holds that $\mathcal{F}_{t+1}:=\mathcal{F}_t$, $\mathcal{A}_{t+1}^1:=\mathcal{A}_{t}^1$, $\mathcal{A}_{t+1}^2:=\mathcal{A}_{t}^2$, $r_{t+1}(F):=r_t(F)$ and $\rho_{t+1}(F):=\rho_t(F)$ for any $F \in \mathcal{F}_{t+1}$ if not specified otherwise.
The update of the remaining sets is automatic according to \Cref{def_sets_phase_2_with_distinction_k3}.

\begin{alg}\label{strat_phase_2_with_distinction_k3}
 In round $t \ge 0$ if $|\mathcal{F}_t|  \geq (1-\alpha)\cdot  \frac{n}{3}$ terminate and, otherwise, let $v_{t+1}$ be the offered vertex.
\begin{enumerate}[label=(\roman*), itemsep=0.7ex]
\item \label{alg2_i} If $v_{t+1} \in  A$ for some $A \in \mathcal{A}^1_t$, choose $w_{t+1} \in V(A)\setminus\{v_{t+1}\}$ arbitrarily such that $v_{t+1}w_{t+1} \notin E(A)$ and thus create a new $A^2$-structure $A':=A \cup \{v_{t+1}w_{t+1}\}$.
We set $\mathcal{A}^1_{t+1}:=\mathcal{A}^1_t \setminus \{A\}$, $\mathcal{A}^2_{t+1}:=\mathcal{A}^2_t \cup \{ A'\}$ and $r_{t+1}(H):=A'$ for all $H \in r_t^{-1}(A)$.
\item \label{alg2_ii} If $v_{t+1}$ is an endpoint of some $A^2$-structure $A = \{ v_{t+1}u_1,u_1u_2\} \in \mathcal{A}^2_t$, we set $w_{t+1}:=u_2$ and thus create a new triangle $F = A\cup \{ v_{t+1}w_{t+1} \}$.
We set $\mathcal{A}^2_{t+1}:=\mathcal{A}^2_t \setminus \{A\}$, $\mathcal{F}_{t+1}:=\mathcal{F}_{t} \cup \{F\}$, and $r_{t+1}(H):=\emptyset$, $\rho_{t+1}(H):=\emptyset$ for all $H \in r_t^{-1}(A) \cup \{F\}$.
\item \label{alg2_iii} If $v_{t+1} \in F$ for some $F \in \mathcal{F}_t$ with $r_t(F) = \rho_t(F) = \emptyset$, u.a.r.~choose $A \in \mathcal{A}^1_t \cup \mathcal{A}^2_t$ and let $w_{t+1}$ be the appropriate vertex of $A$.
We set $r_{t+1}(F):=A$ and $\rho_{t+1}(F) := V(F)$.
\item \label{alg2_iv} If $v_{t+1} \in F$ for some $F \in \mathcal{F}_t$, with $\rho_t(F)=\{v_{t+1}, u_1, u_2\}$ (i.e.\ $|\rho_t(F)|=3$), choose $w_{t+1}$ to be the appropriate vertex of $A=r_t(F)$. 
We set $\rho_{t+1}(F) := \{u_1, u_2\}$.
\item \label{alg2_v} If $v_{t+1} \in F$ for some  $F \in \mathcal{F}_t$, with $\rho_t(F)=\{v_{t+1}, u\}$ (i.e.\ $|\rho_t(F)|=2$), choose $w_{t+1}$ to be the appropriate vertex of $A=r_t(F)$.
We set $\rho_{t+1}(F) = \{u\}$.
\item \label{alg2_vi} If $v_{t+1} \in F$ for some  $F \in \mathcal{F}_t$, with $\rho_t(F)=\{v_{t+1}\}$ (i.e.\ $|\rho_t(F)|=1$), choose $w_{t+1}$ to be the appropriate vertex of $A=r_t(F)$.
Let  $F'$ and $F''$ be the new triangles obtained from $A$ and $F$ by augmentation. Assuming $A\in\mathcal{A}^j_t$, we let $\mathcal{A}^j_{t+1}:=\mathcal{A}^j_{t} \setminus \{ A \}$, $\mathcal{F}_{t+1}:=(\mathcal{F}_{t} \setminus \{ F\}) \cup \{F',F''\}$ and $r_{t+1}(H):=\emptyset$, $\rho_{t+1}(H):=\emptyset$ for all $H \in r_t^{-1}(A) \cup \{F',F''\}$.
\item \label{alg2_vii} Else make an arbitrary move and forget about it.
\end{enumerate}
\end{alg}

Again we claim that w.h.p.~the algorithm terminates.
This time, the only quantities relevant to us are the total number of rounds as well as the number of triangles at the end.

\begin{prp}\label{running_time_phase_2_with_distinction_k3}
 \Cref{strat_phase_2_with_distinction_k3} with $\alpha = \frac{1}{10^{21}}$ and the initial values given by terminal values of \Cref{strat_phase_1_with_r3_k3}  w.h.p.\ terminates after at most $2.251183 \, n$ rounds.
\end{prp}

Crucially, when the algorithm terminates after $t_2$ rounds indeed $\mathcal{F}_{t_2}$ constitutes a family of pairwise vertex-disjoint triangles in $G_{t_2}$.
To see this note that triangles added in rule \ref{alg2_ii} by completion are in $G_{t_2}$ because they were obtained from an $A^2$ structure by adding the missing edge, which itself might have resulted from an $A^1$-structure by adding one more edge.
Moreover, both triangles added in rule \ref{alg2_vi} by augmentation are in $G_{t_2}$ as the algorithm ensures that all edges for the augmentation procedure illustrated in \Cref{fig:augmentation_triangles} are present.
These edges include the edge $a_1a_2$ from the augmentable structure, the edge $f_1f_2$ from the triangle, the edge $f_3a_1$ added in rule \ref{alg2_iii} (or already in rule \ref{alg1_iv} of Phase 1), and the edges $f_3a_2$, $f_1a_3$, and $f_2a_3$ added in rules \ref{alg2_iv}-\ref{alg2_vi} in arbitrary order.

\begin{proof}
First, we have to define the random variables whose development over time we would like to predict. We let $A^1(t):=|A^1_t|$, $A^2(t):=|A^2_t|$, and $F(t):=|F_t|$.
Moreover, as in the proof of \Cref{running_time_phase_1_with_r3_k3} we also track the number of vertices in triangles without responsibility $D(t):=|\cup_{F \in r_t^{-1}(\emptyset)} V(F)|$ and we classify the other vertices in triangles according to the number $i \in \{1,2,3\}$ of suitable offerings their respective triangle still needs as well as by the type of their corresponding augmentable structure $A^j$ with $j \in \{1,2\}$:
\[
R^{A^j,i}(t):=\left|\bigcup_{F \in r_t^{-1}(\mathcal{A}^j_t) \colon |\rho_t(F)|=i} \rho_t(F) \right| \, .
\]

Let $U_{t_1}$, $M_{t_1}$, $A^2_{t_1}$, $F_{t_1}$, $D_{t_1}$, $R^M_{t_1}$, and $R^A_{t_1}$ be the terminal sets of \Cref{strat_phase_1_with_r3_k3}.
We let $A^1(0):=|U_{t_1}|+|M_{t_1}|$, $A^2(0):=|A^2_{t_1}|$, $F(0):=|F_{t_1}|$, $D(0):=|D_{t_1}|$, $R^{A^1,3}(0):=|R^M_{t_1}|$, $R^{A^2,3}(0):=|R^A_{t_1}|$, and $R^{A^1,2}(0)=R^{A^1,1}(0)=R^{A^2,2}(0)=R^{A^2,1}(0)=0$.
W.h.p.~we can assume that the statement of \Cref{running_time_phase_1_with_r3_k3} holds.

We now inspect \Cref{strat_phase_2_with_distinction_k3} to determine the expected differences of the random variables.
The probability that an $A^1$-structure disappears is $A^1(t)/n$ for the chance of it becoming an $A^2$-structure in \ref{alg2_i} and $R^{A^1,1}(t)/n$ for the chance of augmentation in \ref{alg2_vi}.
Similarly the probability for the disappearance of an $A^2$-structure is  $\tfrac{2}{3} \cdot A^2(t)/n + R^{A^2,1}(t)/n$ for \ref{alg2_ii} and \ref{alg2_vi}.
The probability that the number of triangles increases by one is the probability of augmentation plus regular completion, i.e.~$R^{A^1,1}(t)/n + \tfrac{2}{3} A^2(t)/n + R^{A^2,1}(t)/n$ for \ref{alg2_ii} and \ref{alg2_vi}.
With these at hand we already get

\begin{align*}
   \mathbb{E}[A^1(t+1) - A^1(t) \mid H(t)]&= \left(
    \dfrac{A^1(t)}{n} + \dfrac{R^{A^1,1}(t)}{n} \right) \cdot (-3)\\
  \mathbb{E}[A^2(t+1) - A^2(t) \mid H(t)] &= \dfrac{A^1(t)}{n} \cdot 3
   +   \left( \dfrac{\frac{2}{3} \cdot A^2(t)}{n} + \dfrac{R^{A^2,1}(t)}{n} \right) \cdot (-3) \\
  \mathbb{E}[F(t+1) - F(t) \mid H(t)] &= \left( \dfrac{R^{A^1,1}(t)}{n} + \dfrac{\frac{2}{3} \cdot A^2(t)}{n} + \dfrac{R^{A^2,1}(t)}{n} \right)   \cdot 3 \, .
\end{align*}

For $D(t)$ the difficult bit is to determine the number of triangles that lose their responsibility, because either an $A^2$-structure $A$ became a triangle in \ref{alg2_ii} or some augmentable structure $A$ was augmented in \ref{alg2_vi} by some triangle.
In both cases this corresponds to the size of $r_t^{-1}(A)$, but, crucially, it makes a difference if $A$ is selected directly in \ref{alg2_ii} or if it is selected by some triangle $F \in r_t^{-1}(A)$ in \ref{alg2_vi}.
In the former case the expected size of $r_t^{-1}(A)$ is $\tau_2$, where with foresight
\[
\tau_j:=\dfrac{R^{A^j,3}(t) + \frac{3}{2}\cdot R^{A^j,2}(t) + 3 \cdot R^{A^j,1}(t)}{A^j(t)}
\]
for $j=1,2$ is the ratio of responsible triangles for $A^j$-structures to $A^j$-structures.
In the latter case we treat the responsible triangles as balls, augmentable structures as bins, and note that we are in the scenario of \Cref{bib_destruction_lemma}, where one ball (responsible triangle) is chosen u.a.r.~and the corresponding bin (augmentable structure) is destroyed.
By \Cref{bib_destruction_lemma} the expected number of triangles (balls) responsible for the chosen augmentable $A^j$-structure (bin) is given by $1
+ \tau_j + O\left(\tfrac{1}{n}\right)$ .
Accounting for the decrease of $D(t)$ in \ref{alg2_iii} and the additional new triangle created in \ref{alg2_vi} we obtain
\begin{align*}
\mathbb{E}[D(t+1) - D(t) \mid H(t)] = \dfrac{D(t)}{n} \cdot (-3) &+  \dfrac{\frac{2}{3} \cdot A^2(t)}{n}   \cdot 3 \left(1+ \tau_2 \right) \\
&+ \dfrac{R^{A^1,1}(t)}{n}  \cdot 3\left( 2 + \tau_1 + O\left(\dfrac{1}{n}\right) \right) \\
&+  \dfrac{R^{A^2,1}(t)}{n}  \cdot 3 \left( 2 +  \tau_2 +  O\left(\dfrac{1}{n}\right) \right) \, .
\end{align*}
Now we proceed to the six equations $\mathbb{E}[R^{A^j,i}(t+1) - R^{A^j,i}(t) \mid H(t)]$ for $j=1,2$ and $i=1,2,3$.
Naturally, vertices are added to $R^{A^j,3}$ if a triangle without responsibility chooses an $A^j$-structure in \ref{alg2_iii}, $R^{A^j,3}$ decreases by $3$ and $R^{A^j,2}$ increases by $2$ in \ref{alg2_iv}, $R^{A^j,2}$ decreases by $2$ and $R^{A^j,1}$ increases by $1$ in \ref{alg2_v}, and, finally, vertices from $R^{A^j,1}$ perform augmentation in \ref{alg2_vi}.
But also in \ref{alg2_i} responsible triangles are shifted from an $A^1$-structure to an $A^2$-structure and in \ref{alg2_ii} responsibilities become void.

With \ref{alg2_iii} the probability of choosing an $A^j$-structure is
\[
\dfrac{D(t)}{n} \cdot \dfrac{A^j(t)}{A^1(t) + A^2(t)} \, .
\]
The expected number of vertices in triangles that are responsible for an $A^j$-structure and in $\rho_t^{-1}(i)$ is
\[\dfrac{3R^{A^j,i}(t)}{A^j(t)} \, , \]
which is relevant for \ref{alg2_ii},\ref{alg2_iv}, and \ref{alg2_v}.
Lastly, in \ref{alg2_vi}, when an augmentation happens, the expected number of vertices in triangles that are responsible for the augmented $A^j$-structure and in $\rho_t^{-1}(i)$ is again calculated with \Cref{bib_destruction_lemma} and turns out as 
\[1 + \dfrac{3R^{A^j,1}(t)}{A^j(t)} + O\left( \dfrac{1}{n}\right) \, . \]
Summing all these up and using the probabilities we obtained earlier we get
\begin{align*}
    \mathbb{E}[R^{A^1,3}(t+1) - R^{A^1,3}(t) \mid H(t)] &=   \dfrac{D(t)}{n} \cdot \dfrac{A^1(t)}{A^1(t) + A^2(t)} \cdot 3 
   + \dfrac{R^{A^1,3}(t)}{n} \cdot (-3)\\
    &\hspace{1cm} + \left( \dfrac{A^1(t)}{n} + \dfrac{R^{A^1,1}(t)}{n} \right)   \cdot \left( -\dfrac{3R^{A^1,3}(t)}{A^1(t)} \right) \\
    \mathbb{E}[R^{A^1,2}(t+1) - R^{A^1,2}(t) \mid H(t)] &= \dfrac{R^{A^1,3}(t)}{n} \cdot 2
   + \dfrac{R^{A^1,2}(t)}{n} \cdot (-2) \\
    &\hspace{1cm} + \left( \dfrac{A^1(t)}{n} + \dfrac{R^{A^1,1}(t)}{n} \right)   \cdot \left( -  \dfrac{3R^{A^1,2}(t)}{A^1(t)} \right)\\
    \mathbb{E}[R^{A^1,1}(t+1) - R^{A^1,1}(t) \mid H(t)] &=  \dfrac{R^{A^1,2}(t)}{n} \cdot 1
     + \dfrac{A^1(t)}{n} \cdot \left( -\dfrac{3R^{A^1,1}(t)}{A^1(t)} \right)\\
     &\hspace{1cm} +\dfrac{R^{A^1,1}(t)}{n} \cdot \left( -1 - \dfrac{3R^{A^1,3}(t)}{A^1(t)} - O\left( \dfrac{1}{n} \right) \right)
     \\
    \mathbb{E}[R^{A^2,3}(t+1) - R^{A^2,3}(t) \mid H(t)] &=  \dfrac{D(t)}{n} \cdot \dfrac{A^2(t)}{A^1(t) + A^2(t)}  \cdot 3 + \dfrac{R^{A^2,3}(t)}{n} \cdot (-3)\\
     &\hspace{1cm} + \dfrac{A^1(t)}{n}  \cdot   \dfrac{3R^{A^1,3}(t)}{A^1(t)} \\
    &\hspace{1cm} + \left( \dfrac{\frac{2}{3} \cdot A^2(t)}{n} + \dfrac{R^{A^2,1}(t)}{n} \right)   \cdot \left( -  \dfrac{3R^{A^2,3}(t)}{A^2(t)} \right)\\
    \mathbb{E}[R^{A^2,2}(t+1) - R^{A^2,2}(t) \mid H(t)] &= \dfrac{R^{A^2,3}(t)}{n} \cdot 2
   + \dfrac{R^{A^2,2}(t)}{n} \cdot (-2) \\
    &\hspace{1cm} + \dfrac{A^1(t)}{n}  \cdot  \dfrac{3R^{A^1,2}(t)}{A^1(t)}  \\
    &\hspace{1cm} + \left( \dfrac{\frac{2}{3} \cdot A^2(t)}{n} + \dfrac{R^{A^2,1}(t)}{n} \right)   \cdot \left( -  \dfrac{3R^{A^2,2}(t)}{A^2(t)} \right)\\
    \mathbb{E}[R^{A^2,1}(t+1) - R^{A^2,1}(t) \mid H(t)] &=  \dfrac{R^{A^2,2}(t)}{n} \cdot 1
     + \dfrac{A^1(t)}{n}  \cdot    \dfrac{3R^{A^1,1}(t)}{A^1(t)} \\
     &\hspace{1cm} + \dfrac{\frac{2}{3} \cdot A^2(t)}{n}  \cdot \left( -  \dfrac{3R^{A^2,1}(t)}{A^2(t)} \right)\\
      &\hspace{1cm}  + \dfrac{R^{A^2,1}(t)}{n}  \cdot \left( - 1 -  \dfrac{3R^{A^2,3}(t)}{A^2(t)} - O\left( \dfrac{1}{n} \right) \right).
\end{align*}

From the expected difference equations we create a system of differential equations shown in \Cref{fig:system_2}.

\begin{figure}[htbp]
    \centering
\begin{align*}
   a^{1\prime}(x) &= - 3 \cdot a^1(x) - 3 \cdot r^{A^1,1}(x) \\
   a^{2\prime}(x) &= 3 \cdot a^1(x) - 2 \cdot a^2(x) - 3 \cdot r^{A^2,1}(x)\\
   f'(x) &= 3 \cdot r^{A^1,1}(x) + 2 \cdot a^2(x) + 3 \cdot r^{A^2,1}(x)\\
   d'(x) &= -3 \cdot  d(x)  
   + 3 \cdot r^{A^1,1}(x)  \cdot \left(2 +  \dfrac{r^{A^1,3}(x) + \frac{3}{2} \cdot r^{A^1,2}(x) + 3 \cdot r^{A^1,1}(x)}{a^1(x)} \right)\\
   &\hspace{1cm}+ 2 \cdot a^2(x)  \cdot \left(1 +  \dfrac{r^{A^2,3}(x) + \frac{3}{2} \cdot r^{A^2,2}(x) + 3 \cdot r^{A^2,1}(x)}{a^2(x)} \right)\\
   &\hspace{1cm}+ 3 \cdot r^{A^2,1}(x) \cdot \left(2 +  \dfrac{r^{A^2,3}(x) + \frac{3}{2} \cdot r^{A^2,2}(x) + 3 \cdot r^{A^2,1}(x)}{a^2(x)} \right)\\
   r^{A^1,3\prime}(x) &= 3 \cdot \frac{a^1(x)}{a^1(x) + a^2(x)} \cdot  d(x) - 3 \cdot r^{A^1,3}(x) - 3 \cdot \left(a^1(x) + r^{A^1,1}(x) \right) \cdot \frac{r^{A^1,3}(x)}{a^1(x)}\\
   r^{A^1,2\prime}(x) &= 2 \cdot r^{A^1,3}(x) - 2 \cdot r^{A^1,2}(x) - 3 \cdot \left(a^1(x) + r^{A^1,1}(x) \right)\cdot \frac{r^{A^1,2}(x)}{a^1(x)}\\
   r^{A^1,1\prime}(x) &= r^{A^1,2}(x) - 3 \cdot a^1(x)  \cdot \frac{r^{A^1,1}(x)}{a^1(x)} - 3 \cdot r^{A^1,1}(x) \cdot \left( \frac{1}{3} + \frac{r^{A^1,1}(x)}{a^1(x)} \right)\\
   r^{A^2,3\prime}(x) &= 3 \cdot \frac{a^2(x)}{a^1(x) + a^2(x)} \cdot  d(x) - 3 \cdot r^{A^2,3}(x) 
   + 3 \cdot r^{A^1,3}(x)\\
   &\hspace{1cm}- \left(2 \cdot a^2(x) + 3 \cdot r^{A^2,1}(x)\right) \cdot \frac{r^{A^2,3}(x)}{a^2(x)}\\
   r^{A^2,2\prime}(x) &= 2 \cdot r^{A^2,3}(x) - 2 \cdot r^{A^2,2}(x) 
    + 3 \cdot r^{A^1,2}(x) 
   - \left(2 \cdot a^2(x) + 3 \cdot r^{A^2,1}(x)\right)\cdot \frac{r^{A^2,2}(x)}{a^2(x)}\\
   r^{A^2,1 \prime}(x) &= r^{A^2,2}(x)
    + 3 \cdot r^{A^1,1}(x) 
   - 2 \cdot a^2(x)  \cdot \frac{r^{A^2,1}(x)}{a^2(x)}
   -  3 \cdot r^{A^2,1}(x) \cdot \left( \frac{1}{3} + \frac{r^{A^2,1}(x)}{a^2(x)} \right)
\end{align*}

\begin{align*}
 a^1(0) &= \frac{\sqrt{3}}{e^{\sqrt{2} \cdot \arctan(\sqrt{2})}}\\ 
 a^2(0) &= \frac{3 - \sqrt{3}}{e^{\sqrt{2} \cdot \arctan(\sqrt{2})}}\\ 
 f(0) &= 1 - \frac{3}{e^{\sqrt{2} \cdot \arctan(\sqrt{2})}}\\ 
 d(0) &= \frac{1}{4} - \frac{\frac{18}{17}}{\sqrt{3}\cdot e^{\sqrt{2} \cdot \arctan(\sqrt{2})}} 
+  \frac{\frac{39}{68} \cdot \cos( s) + \frac{27\sqrt{15}}{340} \cdot \sin( s)}{e^{\frac{7\sqrt{2}}{4} \cdot \arctan(\sqrt{2})}} \\
r^{A^1,3}(0) &= \frac{3}{8} - \frac{\frac{36}{17}}{\sqrt{3}\cdot e^{\sqrt{2} \cdot \arctan(\sqrt{2})}} 
  - \frac{\frac{99}{136} \cdot \cos(s) - \frac{57\sqrt{15}}{680} \cdot \sin( s)}{e^{\frac{7\sqrt{2}}{4} \cdot \arctan(\sqrt{2})}}\\
r^{A^2,3}(0) &=\frac{3}{8} + \frac{\frac{54}{17} - 3\sqrt{3}}{\sqrt{3}\cdot e^{\sqrt{2} \cdot \arctan(\sqrt{2})}} 
+ \frac{\frac{21}{136} \cdot \cos(s) - \frac{111\sqrt{15}}{680} \cdot \sin( s)}{e^{\frac{7\sqrt{2}}{4} \cdot \arctan(\sqrt{2})}}\\
   r^{A^1,2}(0) &= r^{A^1,1}(0) = r^{A^2,2}(0) = r^{A^2,1}(0) = 0
\end{align*}
    \caption{System of differential equations for \Cref{strat_phase_2_with_distinction_k3} with initial values given by \Cref{running_time_phase_1_with_r3_k3} and $s = \frac{\sqrt{30}}{4}\cdot \arctan(\sqrt{2})$.}
    \label{fig:system_2}
\end{figure}

To apply the DEM we need to make sure that boundedness and trend hypothesis are fulfilled with $\beta(n) = \log(n)^2$ and $\gamma(n) = O\left(\tfrac{1}{n^2}\right)$.
This is essentially immediate for most variables and we only need to show that $D$ does not change too much.
For this we note that in rule \ref{alg2_iii} of \Cref{strat_phase_2_with_distinction_k3} a (modified) Balls into Bins model is hidden, where the triangles act as balls and the augmentable structures act as bins.
The difference is that before the end of Phase~2, all but $\alpha\cdot \frac{n}{3}$ bins will disappear.
So when determining the maximum load of a bin during Phase~2 we can assume that we only have $\alpha\cdot \frac{n}{3}$ bins and that in all $t_2 = x_2 \cdot n$ rounds \ref{alg2_iii} applies.
By \Cref{bib_max_load_lemma} for $b = x_2 \cdot n$, $m = \frac{\alpha}{3} \cdot n$ we get that the probability that a single bin gets at least $\log(n)$ balls is at most $O\left(\frac{1}{n^2}\right)$.

Moreover, keeping in mind the scaling  $x = \frac{t}{n}$ and sufficiently small real numbers $a^{1*}, a^{2*} > 0$ with $a^{1*} + a^{2*} \le \alpha$, the functions on the right hand side of \Cref{fig:system_2} are Lipschitz in the connected open set
\begin{equation*}
\mathcal{D} \coloneqq \underbrace{\left(-\frac{1}{8}, 3\right)}_{\text{for } x} \times 
\underbrace{\left(a^{1*}, 1\right)}_{\text{for } a^1(x)} \times 
 \underbrace{\left(a^{2*}, 1\right)}_{\text{for } a^2(x)} \times 
 \underbrace{(0, 1)}_{\text{for } f(x)} \times 
  \underbrace{(0, 1)}_{\text{for } d(x)} \times 
   \underbrace{\left(-\frac{1}{8}, 1\right)^6}_{\text{for each } r^{A^j,i}(x)}. 
\end{equation*}
Note that we can use the terminal values from \Cref{strat_phase_1_with_r3_k3} with $\lambda(n) \coloneqq n^{-\frac{1}{4}}$ as by \Cref{running_time_phase_1_with_r3_k3} they are concentrated with error $O(n^{3/4})$.

Then we get that all random variables are concentrated around their expected values, which are given by the solution of the ODE-system, with a deviation of at most $O(n^{\frac{3}{4}})$.
Numerical calculations with a precision of 40 digits show that  $f(x) \geq 1 - \frac{1}{2 \cdot 10^{21}}$  before $x$ passes $2.251183$ and that $a^{1*}=a^{2*}=10^{-30}$ are small enough.
So by the DEM we get that w.h.p.
\begin{equation*}
\hspace{2.5cm}  F(2.251183\, n) \geq \left( 1 - \frac{1}{2 \cdot 10^{21}}\right) \cdot n + O(n^{\frac{3}{4}}) \geq \left( 1 - \frac{1}{10^{21}}\right) \cdot n
\end{equation*}
and hence the algorithm terminates.
\end{proof}

\subsection{Phase 3: Building the Final $\alpha\, \frac{n}{3}$ $K_3$'s in $0.002635 \, n$  Rounds}\label{subsection_k3_phase3}

For the final steps of our triangle factor construction, we use an algorithm which constitutes an adaptation of the switching-approach that is presented in Chapters~3.2 and~3.3 of \cite{ben2020very}.
We note upfront that the analysis will not be done by means of the DEM, but concentration inequalities from \Cref{sec:concentration} will be used instead.
The reason for this is that the DEM-based approach from Phase~2 does not work anymore if the number of augmentable structures becomes sublinear (violation of both boundedness and Lipschitz hypothesis).
However, since only a small number of structures is left to cover, a somewhat coarser analysis does not have a dramatic effect on overall running time considerations.

The initial setup of Phase 3 is as follows. For some $\alpha \le \frac{1}{10^{21}}$, we take the set $F$ that consists of $\left( 1 - \alpha \right) \cdot n$ vertices that are contained in complete triangles  and the set $A$ which contains the $\alpha \cdot n$ vertices that are still contained in either $A^1$- or $A^2$-structures.
We note that we will use at most one edge of any augmentable structure and, hence, we treat them as if they all were $A^1$-structures and thus the former distinction does not have to be maintained.

From now on, all new triangles will be created exclusively by augmentation. Moreover, the notion of responsibility will be understood a bit differently now. That is, we divide the set of complete triangles in parts of equal size and decide beforehand which augmentable structure a triangle will take care of once one of its vertices is offered and we refer to this assignment as \emph{responsibility}, despite the potential absence of any edges from this triangle to the augmentable structure.

This clean-up phase will be composed out of many sub-phases, which we will refer to as Phase 3.$S$ for $S \in \mathbb{N}_{\geq 1}$.
The following nomenclature will be used in its analysis:

\begin{defi}
For Phase 3.$S$, with $S \in \mathbb{N}_{\geq 1}$ at the start of Phase 3.$S$ we let $\mathcal{F}_{3.S}$ be the triangles and $\mathcal{A}_{3.S}$ the augmentable structures.
We also set $m_{3.S} \coloneqq \left\lfloor \frac{|\mathcal{F}_{3.S}|}{|\mathcal{A}_{3.S}|}\right\rfloor$ as the number of triangles responsible for an augmentable structure and
$t_{3.S} \coloneqq  \left\lceil n \cdot m_{3.S}^{-\frac{1}{8}} \right\rceil$ for the number of rounds in this phase.
\end{defi}

\begin{figure}[htbp]
    \centering
    \begin{tikzpicture}[scale=0.85]
 \draw[very thick,gray!30, dashed](1.6,1.7) -- (1.6,-2.6);
  \node at (0,1.5) {\scriptsize $a_1^{(1)}$};
  \node at (1,1.5) {\scriptsize $a_2^{(1)}$};
   \node at (0, 0.2) {\scriptsize $a_3^{(1)}$};
    \draw[very thick](0,1) -- (1,1);
\draw[fill=gray!30] (0,1) circle (1.5mm);
\draw[fill=gray!30] (1,1) circle (1.5mm);
\draw[fill=gray!30] (0.5, 0.2) circle (1.5mm);
  \node at (0,-1) {\scriptsize $a_1^{(2)}$};
  \node at (1,-1) {\scriptsize $a_2^{(2)}$};
   \node at (0, -2.3) {\scriptsize $a_3^{(2)}$};
 \draw[very thick](0,-1.5) -- (1,-1.5);
\draw[fill=gray!30] (0,-1.5) circle (1.5mm);
\draw[fill=gray!30] (1,-1.5) circle (1.5mm);
\draw[fill=gray!30] (0.5, -2.3) circle (1.5mm);
\node at (2.4,1.5) {\scriptsize $f_1^{(1,1)}$};
  \node at (3.7,1.5) {\scriptsize $f_2^{(1,1)}$};
   \node at (2.4, 0.2) {\scriptsize $f_3^{(1,1)}$};
    \draw[very thick](2.5,1) -- (3.5,1) -- (3, 0.2) -- cycle;
\draw[fill=gray!30] (2.5,1) circle (1.5mm);
\draw[fill=gray!30] (3.5,1) circle (1.5mm);
\draw[fill=gray!30] (3, 0.2) circle (1.5mm);
  \node at (2.4,-1) {\scriptsize $f_1^{(2,1)}$};
  \node at (3.7,-1) {\scriptsize $f_2^{(2,1)}$};
   \node at (2.4, -2.3) {\scriptsize $f_3^{(2,1)}$};
 \draw[very thick](2.5,-1.5) -- (3.5,-1.5) -- (3, -2.3) -- cycle;
\draw[fill=gray!30] (2.5,-1.5) circle (1.5mm);
\draw[fill=gray!30] (3.5,-1.5) circle (1.5mm);
\draw[fill=gray!30] (3, -2.3) circle (1.5mm);
\node at (4.9,1.5) {\scriptsize $f_1^{(1,2)}$};
  \node at (6.2,1.5) {\scriptsize $f_2^{(1,2)}$};
   \node at (4.9, 0.2) {\scriptsize $f_3^{(1,2)}$};
    \draw[very thick](5,1) -- (6,1) -- (5.5, 0.2) -- cycle;
\draw[fill=gray!30] (5,1) circle (1.5mm);
\draw[fill=gray!30] (6,1) circle (1.5mm);
\draw[fill=gray!30] (5.5, 0.2) circle (1.5mm);
  \node at (4.9,-1) {\scriptsize $f_1^{(2,2)}$};
  \node at (6.2,-1) {\scriptsize $f_2^{(2,2)}$};
   \node at (4.9, -2.3) {\scriptsize $f_3^{(2,2)}$};
 \draw[very thick](5,-1.5) -- (6,-1.5)-- (5.5, -2.3) -- cycle;
\draw[fill=gray!30] (5,-1.5) circle (1.5mm);
\draw[fill=gray!30] (6,-1.5) circle (1.5mm);
\draw[fill=gray!30] (5.5, -2.3) circle (1.5mm);
\node at (7.4,1.5) {\scriptsize $f_1^{(1,3)}$};
  \node at (8.7,1.5) {\scriptsize $f_2^{(1,3)}$};
   \node at (7.4, 0.2) {\scriptsize $f_3^{(1,3)}$};
    \draw[very thick](7.5,1) -- (8.5,1) -- (8, 0.2) -- cycle;
\draw[fill=gray!30] (7.5,1) circle (1.5mm);
\draw[fill=gray!30] (8.5,1) circle (1.5mm);
\draw[fill=gray!30] (8, 0.2) circle (1.5mm);
  \node at (7.4,-1) {\scriptsize $f_1^{(2,3)}$};
  \node at (8.7,-1) {\scriptsize $f_2^{(2,3)}$};
   \node at (7.4, -2.3) {\scriptsize $f_3^{(2,3)}$};
 \draw[very thick](7.5,-1.5) -- (8.5,-1.5)-- (8, -2.3) -- cycle;
\draw[fill=gray!30] (7.5,-1.5) circle (1.5mm);
\draw[fill=gray!30] (8.5,-1.5) circle (1.5mm);
\draw[fill=gray!30] (8, -2.3) circle (1.5mm);
  \node at (11.5,0.25) {\scriptsize (no responsibility)};
    \draw[very thick](11,-0.25) -- (12,-0.25) -- (11.5, -1.05) -- cycle;
\draw[fill=gray!30] (11,-0.25) circle (1.5mm);
\draw[fill=gray!30] (12,-0.25) circle (1.5mm);
\draw[fill=gray!30] (11.5, -1.05) circle (1.5mm);
  \end{tikzpicture}
    \caption{Example with $|\mathcal{A}_{3.S}| = 2$, $|\mathcal{F}_{3.S}| = 7$, which implies that $m_{3.S} = 3$.}
    \label{fig:Algo_3}
\end{figure}
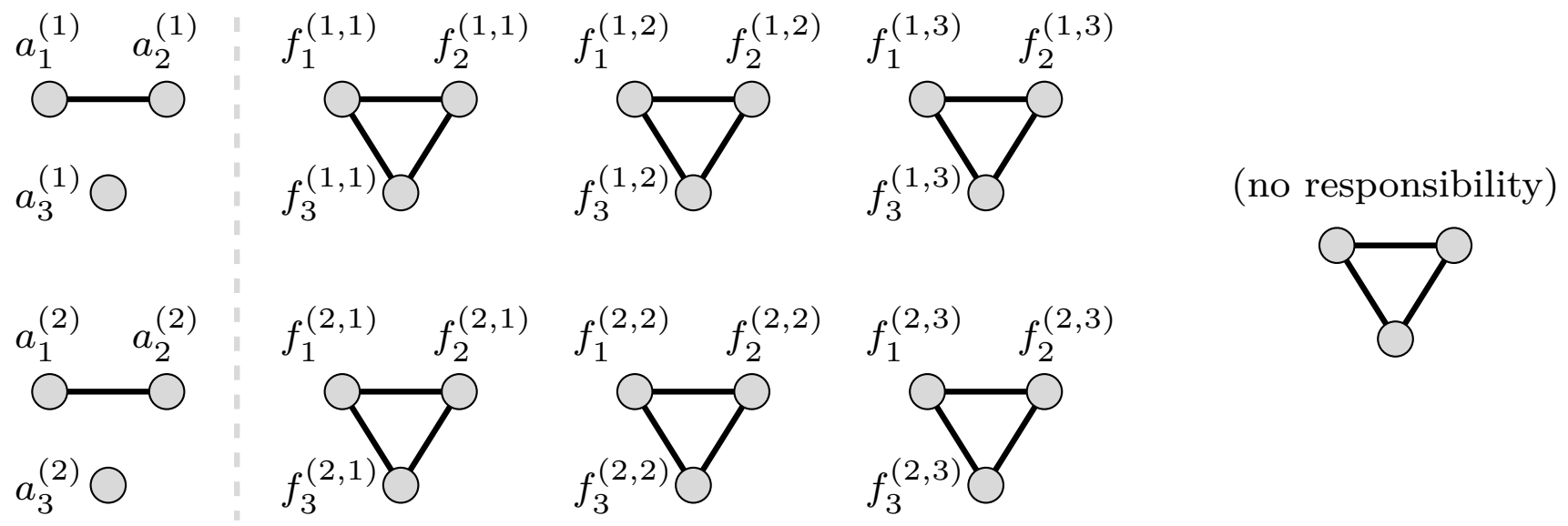

As in Phase~2, during the $S$-th sub-phase of Phase~3, augmentation will be possible once a triangle that is responsible for an augmentable structure gets all of its $4$ suitable offerings. 
The first vertex of a triangle that is offered during Phase~3.$S$ gets labeled $f_3$ and afterwards each vertex of that triangle needs one more offering. 
We illustrate the situation graphically below in \Cref{fig:Algo_3}.
For the formal analysis of Phase~3.$S$, we are interested in the event that augmentation can happen for specifically the $i$-th augmentable structure and its $j$-th responsible triangle $T^{(i,j)}$. 
\begin{defi}
For Phase 3.$S$, with $S \in \mathbb{N}_{\geq 1}$, for $i \in [|\mathcal{A}_{3.S}|]$, $j \in [m_{3.S}]$, and $t_{3.S} \in \mathbb{N}$, define the event
\begin{align*}
   \mathcal{E}_S^{(i,j)} &\coloneqq \text{in $t_{3.S}$ rounds, the vertices
   of $T^{(i,j)}$ have received all $4$ suitable offerings.}
\end{align*}
\end{defi}

Note that this event does not exclude that the vertices involved receive more than the aforementioned $4$ offerings.
The following algorithm now describes how Builder creates new triangles during Phase 3.
 \begin{alg}\label{strat_phase_3_k3}
For $S \in \mathbb{N}_{\geq 1}$, if $\mathcal{A}_{3.S} = \emptyset$ terminate and, otherwise, do the following in Phase 3.$S$:
\begin{itemize}
    \item Arbitrarily label the augmentable structures in $\mathcal{A}_{3.S}$ as $A^{(i)}$ for $i \in [|\mathcal{A}_{3.S}|]$ and label $|\mathcal{A}_{3.S}| \cdot m_{3.S}$ of the triangles in $\mathcal{F}_{3.S}$ as $T^{(i,j)}$ with $i \in [|\mathcal{A}_{3.S}|]$ and $j \in [m_{3.S}]$. (Call the unlabeled triangles $\mathcal{U}_{3.S} \coloneqq \mathcal{F}_{3.S}\setminus \{T^{(i,j)}: i\in[|\mathcal{A}_{3.S}|], j\in[m_{3.S}]\}$.)
    \item For the next $t_{3.S}$ rounds of the SRGP if $v_{t+1}$ is a vertex of triangle $T^{(i,j)}$ and it did not have all its suitable offerings yet, choose $w_{t+1}$ to be the appropriate vertex of the $i$-th augmentable structure $A^{(i)}$. 
    Otherwise, make an arbitrary move and forget about it.
\item Let $I \subseteq [|\mathcal{A}_{3.S}|]$, be those $i$ for which there is $j \in [m_{3.S}]$ such that $\mathcal{E}_S^{(i,j)}$ has happened.
For each $i \in I$ we pick one such $j_i$ and let $F_i$ and $F_i'$ be the triangles obtained from $A^{(i)}$ and $T^{(i,j_i)}$ by augmentation.
The remaining augmentable structures are $\mathcal{A}_{3.S+1} \coloneqq \{ A^{(i)} \colon i \not\in I \}$ and the triangles are
\begin{align*}
\mathcal{F}_{3.S+1} \coloneqq &\{ T^{(i,j)} \in \mathcal{F}_{3.S} \colon i \not\in I \} \cup \mathcal{U}_{3.S}\\
 &\cup \{ F_i,F_i' \colon i \in I \} \cup \{ T^{(i,j)} \in \mathcal{F}_{3.S} \colon i \in I, j \in [m_{3.S}]  \setminus \{j_i\} \}   \, .
\end{align*}
\end{itemize}
\end{alg}

Our first goal is to lower bound the probability of the event $\mathcal{E}_S^{(i,j)}$ that triangle $T^{(i,j)}$ can augment $A^{(i)}$.

\begin{lm}\label{q_lemma_k3}
For $S \in \mathbb{N}_{\geq 1}$, $\alpha \le 10^{-15}$, $i \in [|\mathcal{A}_{3.S}|]$, and $j \in [m_{3.S}]$
\begin{equation*}
    \mathbb{P}(\mathcal{E}_S^{(i,j)}) \ge  q_{3.S} := 0.4998\, e^{-3} \cdot m_{3.S}^{-1/2} \, .
\end{equation*}
\end{lm}

\begin{proof}
First it is worth noting that $\frac{1}{\alpha}-2 \leq m_{3.S}$.
This is the case because 
\begin{equation*}
m_{3.S} = \left\lfloor \dfrac{|\mathcal{F}_{3.S}|}{|\mathcal{A}_{3.S}|}\right\rfloor 
 \geq ... 
  \geq \left\lfloor \dfrac{|\mathcal{F}_{3.1}|}{|\mathcal{A}_{3.1}|}\right\rfloor 
\ge \left\lfloor \dfrac{\left( 1 - \alpha \right)\cdot n}{\alpha \cdot n}\right\rfloor = \left\lfloor  \frac{1}{\alpha} - 1 \right\rfloor \geq  \frac{1}{\alpha} - 2.
\end{equation*}
Moreover, using that $m_{3.S} \leq n$ and $\alpha \leq \frac{1}{10^9}$, we get
\begin{equation*}
 \left\lceil n^{\frac{7}{8}}\right\rceil 
 =\left\lceil n \cdot n^{-\frac{1}{8}}\right\rceil  
 \leq \left\lceil n \cdot m_{3.S}^{-\frac{1}{8}} \right\rceil = t_{3.S}
   \leq  \left\lceil n \cdot \left( \frac{1}{\alpha}-2\right)^{-\frac{1}{8}} \right\rceil  
   \leq \left\lceil n \cdot 0.1 \right\rceil.
\end{equation*}
Hence for $n$ large enough, we have in particular that  $n \ge t_{3.S}\ge 30000$.

Now we observe that $\mathbb{P}\left[ \mathcal{E}^{(i,j)}_S \right]$
 is at least the probability that in $t_{3.S}$ rounds, the vertices
   of $T^{(i,j)}$ have received \emph{exactly} their $4$ suitable offerings,
and note that there are exactly $(3 \cdot 2) \cdot 3 = 18$ orderings in which the $4$ vertex offerings can happen, since the first vertex offered takes the role of  $f_3$. 

Then we can estimate
\begin{align*}
 \mathbb{P}\left[ \mathcal{E}^{(i,j)}_S \right] 
    &\geq \binom{t_{3.S}}{4} \cdot 18
    \cdot \left( \frac{1}{n} \right)^4  \cdot \left( 1 - \frac{3}{n} \right)^{t_{3.S}-4} \\
     &\geq \dfrac{t_{3.S} \cdot (t_{3.S}-1) \cdot (t_{3.S}-2) \cdot (t_{3.S}-3)}{24} \cdot 18 \cdot \left( \frac{1}{n} \right)^4 \cdot e^{- \frac{\frac{3}{n}}{1 - \frac{3}{n}} \cdot (t_{3.S}-4)} \\
     &\geq \dfrac{0.9999^4 \cdot t_{3.S}^4}{2 \cdot n^4}  \cdot e^{- \frac{3 \cdot (t_{3.S}-4)}{n-3}}\\
       &\geq  0.4998 \cdot m_{3.S}^{-\frac{1}{2}} \cdot e^{- 3}=q_{3.S}.
\end{align*}
Hence $q_{3.S}$ is a suitable lower bound for $\mathbb{P}\left[ \mathcal{E}^{(i,j)}_S \right]$.
\end{proof}

Next, we address the question of how much the set of vertices in augmentable structures actually shrinks during a single sub-phase. An expression that will be useful in many of the following calculations is 
\begin{equation*}
   m_{3.S} \cdot q_{3.S} =m_{3.S} \cdot 0.4998\, e^{-3} \cdot m_{3.S}^{-1/2}
   = \frac{0.4998}{e^3} \cdot \sqrt{ \left\lfloor\frac{|\mathcal{F}_{3.S}|}{|\mathcal{A}_{3.S}|}\right\rfloor}
    \geq \frac{0.4998}{e^3} \cdot \sqrt{\frac{1}{\alpha} - 2}.
\end{equation*}
In particular, for $\alpha \leq 10^{-15}$, it holds that $ m_{3.S} \cdot q_{3.S} \geq 786887$.

Now we state and prove a lemma that investigates the development of the set of augmentable structures depending on the aforementioned product.
\begin{lm}\label{technical_lemma_k3}
    For $\alpha \le \frac{1}{10^{15}}$ and $S \in \mathbb{N}_{\geq 1}$ after Phase 3.$S$
\begin{itemize}
    \item if $\frac{m_{3.S} \cdot q_{3.S}}{128} \geq \log(n)$ we have $|\mathcal{A}_{3.(S+1)}|=0$ w.h.p.~and
    \item if $\frac{m_{3.S} \cdot q_{3.S}}{128} < \log(n)$ we have $|\mathcal{A}_{3.(S+1)}|\le 2 \, e^{-\frac{m_{3.S} \cdot q_{3.S}}{512}} \cdot n$ with probability at least $1-e^{-\sqrt{n}}$.
\end{itemize}
\end{lm}

\begin{proof}
Since the statement is intended to apply to all Phases 3.$S$, we will suppress most of the indices in the following calculations.
For $i \in [|\mathcal{A}|]$, we define the random variable 
$X_i \coloneqq \left| \left\{ j \in [m] \colon \mathds{1}_{\mathcal{E}^{(i,j)}}=1 \right\} \right|$.
It follows with \Cref{q_lemma_k3} that
\begin{equation*}
  \mathbb{E}[X_i] = \sum_{j \in [m]} \mathbb{E}\left[ \mathds{1}_{\mathcal{E}^{(i,j)}}  \right] = \sum_{j \in [m]} \mathbb{P}\left[ \mathcal{E}^{(i,j)}  \right] \geq m \cdot q.
\end{equation*}
To gain further insights regarding the concentration of $X_i$, we use Talagrand's inequality (cf.\ \Cref{lem:talagrand}) with $c = 1$, $g = 4$, and $\lambda = \frac{\mathbb{E}[X_i]}{2}$ to estimate
\begin{align*}
    \mathbb{P}\left[ X_i < \frac{\mathbb{E}[X_i]}{2} - 60 \cdot \sqrt{4 \cdot \mathbb{E}[X_i]} \right] 
    &\leq \mathbb{P} \left[ |X_i - \mathbb{E}[X_i]| > \frac{\mathbb{E}[X_i]}{2}+ 60 \cdot \sqrt{4 \cdot \mathbb{E}[X_i]} \right]\\
    &\leq 4 \cdot e^{- \frac{\lambda^2}{8 \cdot g \cdot \mathbb{E}[X_i]}} = 4 \cdot e^{- \frac{\mathbb{E}[X_i]}{128}} \leq 4 \cdot e^{- \frac{m \cdot q}{128}}.
\end{align*}
For this bound to be informative, it must hold that
\begin{equation*}
    \frac{\mathbb{E}[X_i]}{2} - 60 \cdot \sqrt{4 \cdot \mathbb{E}[X_i]} = \sqrt{\mathbb{E}[X_i]} \cdot \left( \frac{\sqrt{\mathbb{E}[X_i]}}{2} - 120\right) \geq \sqrt{\mathbb{E}[X_i]} \cdot \left( \frac{\sqrt{m \cdot q}}{2} - 120\right) > 0,
\end{equation*}
which is indeed the case as $\mathbb{E}[X_i] \geq m \cdot q \geq 786887$ and thus $ \left( \frac{\sqrt{m \cdot q}}{2} - 120\right) > 300$.

Now we define the random variable 
$Y \coloneqq | \{ i \in [|\mathcal{A}|] \colon X_i=0\} | \, $ and derive 
\begin{equation*}
    \mathbb{E}[Y] = 
    \sum_{i \in [|\mathcal{A}|]} \mathbb{P}\left[ X_i = 0 \right] 
    \leq  \sum_{i \in [|\mathcal{A}|]} \mathbb{P}\left[ X_i <\frac{\mathbb{E}[X_i]}{2} - 60 \cdot \sqrt{4 \cdot \mathbb{E}[X_i]}\right] \leq |\mathcal{A}| \cdot 4 \cdot e^{- \frac{m \cdot q}{128}}.
\end{equation*}
We conduct a case distinction based on the size of the exponent in the previous expression.

\textbf{Case 1}: $\frac{m \cdot q}{128} \geq \log(n)$.
As $m \ge \log(n)$ we have $|\mathcal{A}| \le \frac{n}{\log(n)}$ and
\begin{equation*}
    \mathbb{E}[Y] \leq  \frac{n}{\log(n)} \cdot 4 \cdot e^{-\log(n)} = 
    \frac{n}{\log(n)} \cdot 4 \cdot \frac{1}{n} = o(1).
\end{equation*}
By Markov's inequality $\mathbb{P} [Y \geq 1] = o(1)$ and w.h.p.\  the $K_3$-factor is completed.

\textbf{Case 2}: $\frac{m \cdot q}{128} < \log(n)$.
Using Azuma's inequality (cf.\ \Cref{lem:azuma}), with $c = 1$, $\lambda =  e^{- \frac{1}{4} \cdot \frac{m \cdot q}{128}} \cdot n > e^{- \frac{1}{4} \cdot \log(n)} \cdot n = n^{\frac{3}{4}}$, and $t_{3.S} \le  n$, we get
\begin{equation*}
    \mathbb{P}\left[ Y \geq \mathbb{E}[Y]  + \lambda \right] \leq e^{- \frac{\lambda^2}{2 \cdot t_{3.S}}} \leq e^{- \frac{n^{3/2}}{2 \cdot  n}} \leq e^{- \sqrt{n}}.
\end{equation*}
Hence with probability at least $1 - e^{-\sqrt{n}}$ we have $Y \leq \mathbb{E}[Y]  + \lambda$ and thus
\begin{align*}
   |\mathcal{A}_{3.(S+1)}| \leq \mathbb{E}[Y]  + \lambda
     \leq \tfrac{\alpha \cdot n}{3} \cdot 4 \cdot e^{- \frac{m \cdot q}{128}} + e^{- \frac{1}{4} \cdot \frac{m \cdot q}{128}} \cdot n \leq 2 \, e^{- \frac{m \cdot q}{512}} \cdot n \, ,
\end{align*}
where the last inequality uses the bound on $\alpha$.
\end{proof}

 Based on the previous result, the next lemma shows that $\frac{3|\mathcal{A}_{3.S}|}{n}$ shrinks (at least) exponentially with $S$ and becomes zero within $O(\log(\log(n)))$ sub-phases.

\begin{lm}\label{induction_lemma_k3}
For $\alpha \le \frac{1}{10^{15}}$ it holds w.h.p.~that for all $S \in \mathbb{N}_{\geq 1}$ with $\frac{m_{3.S} \cdot q_{3.S}}{128} < \log(n)$ we have
\begin{equation*}
\frac{3|\mathcal{A}_{3.S}|}{n} \leq \alpha^{S}
\end{equation*}
and the whole Phase 3 consists of at most $\log(\log(n))$ sub-phases and succeeds.
\end{lm}

\begin{proof} 
We proceed by induction on $S$.
The base case holds as $3|\mathcal{A}_{3.1}| = \alpha \cdot n$.
For the inductive step we use $3|\mathcal{F}_{3.S}| \ge (1- \alpha)\cdot n$ to get
\begin{equation*}
   m_{3.S} \cdot q_{3.S}
     =  0.4998\, e^{-3} \cdot \sqrt{ \left\lfloor \frac{|\mathcal{F}_{3.S}|}{|\mathcal{A}_{3.S}|}\right\rfloor}
    \geq 0.4998\, e^{-3} \cdot \sqrt{ \frac{(1-2\alpha) \cdot n}{3|\mathcal{A}_{3.S}|}}
    \geq  0.4997\, e^{-3} \cdot \sqrt{\frac{n}{3|\mathcal{A}_{3.S}|}}.
\end{equation*}
Assuming $\frac{m_{3.S} \cdot q_{3.S}}{128} < \log(n)$ and as $\frac{m_{3.S} \cdot q_{3.S}}{128} = \frac{0.4998}{e^3 \cdot 128} \cdot \sqrt{m_{3.S}}$ is monotonically increasing in $S$, we can apply the induction hypothesis to estimate
\begin{equation*}
    \sqrt{\frac{n}{3|\mathcal{A}_{3.S}|}} \geq \sqrt{\left(\frac{1}{\alpha}\right)^S} = \alpha^{-\frac{S}{2}}.
\end{equation*}
With this, \Cref{technical_lemma_k3}, and the estimate above for $m_{3.S} \cdot q_{3.S}$ with probability at least $1-e^{-\sqrt{n}}$ we get
\begin{equation*}
    \frac{3|\mathcal{A}_{3.(S+1)}|}{n} \leq  6\, e^{-\frac{m_{3.S} \cdot q_{3.S}}{512}} \leq 6 \, e^{ -\frac{0.4998}{e^3 \cdot 512} \cdot \sqrt{\frac{n}{3|\mathcal{A}_{3.S}|}}}  \le 
    e^{ \log(6) -\frac{0.4998}{e^3 \cdot 512}  \cdot \alpha^{-\frac{S}{2}}} \, .
\end{equation*}
However, this is at most $ \alpha^{S+1}$, when $S \geq 1$ and $\alpha \le \frac{1}{10^{15}}$.
Therefore, the inductive step fails with probability at most $e^{-\sqrt{n}}$.

From \Cref{technical_lemma_k3} we know that the triangle factor is w.h.p.\ complete, once $\frac{m_{3.S} \cdot q_{3.S}}{128}$ reaches $\log(n)$. 
From above we get
\begin{equation*}
   \frac{m_{3.S} \cdot q_{3.S}}{128}
    \geq  \frac{0.4997}{e^3 \cdot 128} \cdot \sqrt{\frac{n}{3|\mathcal{A}_{3.S}|}}
    \geq  \frac{1}{10^5}\cdot \alpha^{-\frac{S}{2}} \, ,
\end{equation*}
which, for $\alpha \le \frac{1}{10^{15}}$, $S \ge \log(\log(n))$, and $n$ sufficiently large is at least $\log (n)$.

Each inductive step, with the exception of the last one, fails with probability at most $e^{- \sqrt{n}}$ and the last one, by \Cref{technical_lemma_k3}, fails with probability at most $o(1)$.
Hence the  total probability of failure is bounded by
\begin{equation*}
  \hspace{1.85cm}  \log(\log(n)) \cdot e^{-\sqrt{n}} + o(1) = e^{\log\left(\log(\log(n))  \right)- \sqrt{n}} + o(1) = 
    o(1). \qedhere
\end{equation*}
\end{proof}

With this we can now estimate the total number of rounds for Phase 3.

\begin{prp}\label{phase_3_rounds_lemma_k3}
For $\alpha \le \frac{1}{10^{15}}$, \Cref{strat_phase_3_k3} w.h.p.~terminates after at most $\tfrac{10}{9} \, \alpha^{1/8} n$ rounds.
\end{prp}

\begin{proof}
Knowing from \Cref{induction_lemma_k3} that \Cref{strat_phase_3_k3} w.h.p.\ succeeds, it remains to estimate its total duration.
Let $S^* \le \log (\log(n))$ be the smallest $S \in \mathbb{N}$ such that $\frac{m_{3.S} \cdot q_{3.S}}{128} \geq \log(n)$. So by definition of $t_{3.S}$, it holds for $t_3$, the total number of rounds of Phase 3, that
\begin{equation*}
  t_3 = \sum_{S=1}^{S^*}  t_{3.S} = \sum_{S=1}^{S^*}   \left\lceil n \cdot m_{3.S}^{-\frac{1}{8}} \right\rceil \leq  \sum_{S=1}^{S^*}  \left( n \cdot \left\lfloor \frac{|\mathcal{F}_{3.S}|}{|\mathcal{A}_{3.S}|}\right\rfloor^{-\frac{1}{8}} + 1 \right) \leq  \sum_{S=1}^{S^*} n \cdot \left(  \frac{|\mathcal{A}_{3.S}| \cdot 2}{|\mathcal{F}_{3.S}|} \right)^{\frac{1}{8}} +  S^*.
\end{equation*}
By \Cref{induction_lemma_k3} w.h.p.~$3|\mathcal{A}_{3.S}| \leq \alpha^S \cdot n$ and as $3|\mathcal{F}_{3.S}| =  \left( 1 - \alpha\right) \cdot n$ we further estimate with $\alpha \le \frac{1}{10^{15}}$
\begin{align*}
  t_3 &\leq  \sum_{S=1}^{S^*} n \cdot \left(  \frac{\alpha^S \cdot 2}{1 - \alpha} \right)^{\frac{1}{8}} +  S^* 
  \leq n \cdot  \left(  \frac{ 2}{1 - \alpha} \right)^{\frac{1}{8}} \cdot \sum_{S=1}^{\infty}\left( \alpha^{\frac{1}{8}} \right)^S +  S^* \\
  &= n \cdot \left(  \frac{ 2}{1 - \alpha} \right)^{\frac{1}{8}} \cdot \frac{\alpha^{\frac{1}{8}} }{1 - \alpha^{\frac{1}{8}} } +  S^* \leq \tfrac{10}{9} \, \alpha^{1/8} \cdot n \, .
\end{align*}
\end{proof}

\begin{proof}[Proof of \Cref{main_theorem_k3_factor_upper}]
According to \Cref{running_time_phase_1_with_r3_k3}, \Cref{strat_phase_1_with_r3_k3} w.h.p.\ terminates after at most $0.675511\, n$ rounds.
Now we use the results of this algorithm as initial values for Phase 2 and set $\alpha \coloneqq \frac{1}{10^{21}}$.
With this choice of $\alpha$, by \Cref{running_time_phase_2_with_distinction_k3}, \Cref{strat_phase_2_with_distinction_k3} w.h.p.\ terminates after at most $ 2.251183 \, n$ rounds.
Finally, by \Cref{phase_3_rounds_lemma_k3}, using the same value for $\alpha$, Phase 3 consists of at most $10/9 \hspace{0.1cm} \alpha^{1/8} n\le 0.002635 \, n$ rounds.
Hence the total number of rounds it takes to construct a $K_3$-factor is w.h.p.\ at most $2.929329\, n$. 
\end{proof}

\section{Conclusion}\label{section_conclusion}

 As for the upper bounds on $\tau_{\mathcal{K}_3^f}$ and $\tau_{\mathcal{K}_4^f}$, our results are certainly not tight.
 While our three phase distinction significantly simplifies the analysis it is inherently suboptimal to only allow the combination of completion and augmentation in the second phase.
 In \cite[Section 4.3.4.1]{eva_thesis} there is a discussion of some possible improvements.

A search for better lower bounds on $\tau_{\mathcal{K}_k^f}$ would be an intriguing endeavor as well.
Recall that for $k=3$, the proof for $ 1.26575 < \tau_{\mathcal{HAM}}$  also works for $\tau_{\mathcal{K}_3^f}$, because it only uses that the graph is $2$-regular.
So on the one hand, it would be interesting to come up with a proof that is more tailored to triangle factors.
\footnote{A necessary condition for a triangle factor is that every vertex that was never offered has at least two incoming edges. Since the triangles in a triangle factor are pairwise vertex-disjoint, the set of vertices sending edges to unoffered vertices must therefore have size at least twice the number of unoffered vertices.
An exploratory exchange with ChatGPT suggested that this observation yields a lower bound of $\tau_{K_3^f}(n)\geq 1.3426$.
We do not pursue that bound here, as we expect that a stronger one can be obtained by exploiting further structural constraints, for example, that the neighbor of a neighbor of a neighbor of $v$ must be $v$ itself.}
In that regard, it would be particularly gratifying to show that $1.81696 < \tau_{\mathcal{K}_3^f}$, the latter number being the currently best upper bound on $\tau_{\mathcal{HAM}}$.
This would establish that $\tau_{\mathcal{HAM}} < \tau_{\mathcal{K}_3^f}$.
We indeed conjecture the latter to be the case, as a $K_3$-factor is notoriously ``harder to get'' than the equally $2$-regular Hamilton cycle, as can be seen by the minimum degree conditions of  $2n/3$~\cite{hajnalszemeredi} vs. $n/2$~\cite{diracham} and the random graph thresholds $\frac{\log(n)^{1/3}}{n^{2/3}}$~\cite{johansson2008factorsrandomgraphs} vs.  $\frac{\log(n)}{n}$~\cite{posaham}.
    
On the other hand, it would also be interesting to generalize the generic lower bound proof for $2$-regular graphs to $r$-regular graphs with arbitrary $r \geq 3$.
This would definitely constitute a tremendous improvement on the currently best known lower bound on $\tau_{\mathcal{K}_{r+1}^f}$, namely the offline bound.
With $r$ tending to infinity the trivial lower bound $rn/2$ is asymptotically optimal \cite{anastos2025optimally,ben2020very}.

\section*{Declaration of AI use}
AI tools were used in a supporting capacity to create \Cref{fig:plot}, cross-check citations, improve language, and proofread the manuscript. They were not used to develop the mathematical ideas or arguments presented in this work. The authors take full responsibility for the content of this paper.

\bibliographystyle{amsplain}
\bibliography{references}

\end{document}